\documentclass[a4paper]{amsart}
\usepackage{times,mathrsfs,stmaryrd}
\usepackage{pstricks}
\usepackage{amsmath,amssymb}
\usepackage{booktabs} 
\newdimen\hoogte    \hoogte=10pt    
\newdimen\breedte   \breedte=13pt   
\newdimen\dikte     \dikte=0.6pt    
\newenvironment{Young}{\begingroup
       \def\vr{\vrule height0.8\hoogte width\dikte depth 0.2\hoogte}
       \def\fbox##1{\vbox{\offinterlineskip
                    \hrule height\dikte
                    \hbox to \breedte{\vr\hfill$##1$\hfill\vr}
                    \hrule height\dikte}}
       \vtop\bgroup \offinterlineskip \tabskip=-\dikte \lineskip=-\dikte
            \halign\bgroup &\fbox{##\unskip}\unskip  \crcr }
       {\egroup\egroup\endgroup}
\def\Tab(#1){\text{\tiny\begin{Young}#1\cr\end{Young}}}

\usepackage[backref=page]{hyperref}
\renewcommand*{\backref}[1]{}
\renewcommand*{\backrefalt}[4]{{[\tiny%
  \ifcase #1 Not cited.\relax\or Page~#2.%
  \else Pages #2.\fi]}}

\let\len=\ell
\let\<=\langle
\let\>=\rangle
\def\({\big(}
\def\){\big)}
\let\bar=\overline

\def\bijection{\overset{\simeq}{\longrightarrow}}

\let\concat=\sqcup
\let\phi=\varphi
\let\epsilon=\varepsilon

\newcommand{\field}{\Bbbk}

\newcommand{\Comp}[1][n]{\Lambda_{#1}}
\newcommand\SComp[1][n]{\Lambda^\pm_{#1}}
\def\SPart{\Lambda^\oplus_n}
\def\C{\mathbb C}
\def\D{\mathscr D}
\def\E{\mathscr E}

\def\alpp{\alpha^+}
\def\lamp{\lambda^+}
\def\mup{{\mu^+}}
\def\nup{{\nu^+}}
\def\sigmap{\sigma^+}
\def\G{G_{r,n}}
\def\ColG{\mathscr G(r)}

\def\M{\mathcal N}
\def\N{\mathbb N}
\def\P{\mathbf P}
\let\Sect\S
\def\S{\mathscr S}
\def\Words{\mathcal G_{r,n}}

\def\c{\mathbf c}

\def\n{\mathbf{n}}
\def\rn{\mathbf{n}_\zeta}

\def\Sym{\mathfrak S}
\def\s{\mathfrak s}
\def\t{\mathfrak t}
\def\rssTab{\mathcal{T}}
\def\srssTab{\mathcal{T}_\zeta}
\def\T{\mathsf T}

\def\Z{\mathbb Z}
\def\map#1#2{\,{:}\,#1\!\longrightarrow\!#2}

\DeclareMathOperator{\sequiv}{\sim_{\text{sgn}}}
\DeclareMathOperator{\row}{\mathsf{row}}
\DeclareMathOperator{\col}{\mathsf{col}}
\DeclareMathOperator{\comp}{\mathsf{comp}}
\def\mat#1{\big(\begin{smallmatrix}#1\end{smallmatrix}\big)}
\DeclareMathOperator{\Mat}{\mathsf{Mat}}
\DeclareMathOperator{\wt}{\mathsf{wt}}
\DeclareMathOperator{\id}{\mathsf{id}}
\DeclareMathOperator{\Col}{\mathsf{Col}}
\DeclareMathOperator{\Des}{\mathsf{Des}}
\DeclareMathOperator{\Sol}{\mathsf{Sol}}
\DeclareMathOperator{\Rad}{Rad}

\let\emph=\textbf

\numberwithin{equation}{section}
\swapnumbers
\newtheorem{Defn}[equation]{Definition}
\newtheorem{Thm}[equation]{Theorem}
\newtheorem{Prop}[equation]{Proposition}
\newtheorem{Lemma}[equation]{Lemma}
\newtheorem{Cor}[equation]{Corollary}

\newtheorem{Point}[equation]{}

\theoremstyle{remark}
\newtheorem{Remark}[equation]{Remark}

\newenvironment{Example}%
  {\refstepcounter{equation}\trivlist
   \item[\hskip\labelsep\theequation.~\textbf{Example}\space]
   \ignorespaces
  }{\unskip\nobreak\hfil%
    \penalty50\hskip2em\hbox{}\nobreak\hfil$\Diamond$%
    \parfillskip=0pt\finalhyphendemerits=0\penalty-100\endtrivlist}

\makeatletter
{\catcode`\|=\active
  \gdef\set#1{\mathinner{\lbrace\,{\mathcode`\|"8000%
                                   \let|\midvert #1}\,\rbrace}}
}
{\catcode`\|=\active
  \gdef\Set#1{\mathinner{\Big\lbrace\,{\mathcode`\|"8000%
                                   \let|\Midvert #1}\,\Big\rbrace}}
}
\def\midvert{\egroup:\bgroup}
\def\Midvert{\egroup\,\Big|\,\bgroup}

\def\SET[#1]#2|#3|{\Bigg\lbrace\,#2\,\Bigg|\,%
           \vcenter{\hsize #1mm\centering #3}\Bigg\rbrace}

\@tempcnta\time\@tempcntb\time
\divide\@tempcnta by 60
\multiply\count255 by -60 
\advance\@tempcntb by \count255
\edef\Time{\ifnum\@tempcntb<10 0\fi\number\@tempcntb}
\ifnum\@tempcnta>12%
  \advance\@tempcnta by -12
  \edef\Time{\number\@tempcnta:\Time\space PM}%
\else\ifnum\@tempcnta=12
  \edef\Time{\number\@tempcnta:\Time\space PM}%
\else
  \edef\Time{\number\@tempcnta:\Time\space AM}%
\fi\fi

\begin{document}\bibliographystyle{andrew}
\let\@evenfoot\@oddfoot
\makeatother

\title{Cyclotomic Solomon algebras}
\author{Andrew Mathas}
\address{School of Mathematics and Statistics, %
         University of Sydney, NSW 2006, Australia.}
\email{a.mathas@usyd.edu.au}

\author{Rosa C. Orellana}
\address{Department of Mathematics, Dartmouth College, %
         Hanover, NH 03755--3551, USA.}
\email{rosa.c.orellana@dartmouth.edu}


\begin{abstract}
   This paper introduces an analogue of the Solomon descent algebra for the
   complex reflection groups of type $G(r,1,n)$. As with the Solomon
   descent algebra, our algebra has a basis given by sums of
   `distinguished' coset representatives for certain `reflection
   subgroups'. We explicitly describe the structure constants with respect
   to this basis and show that they are polynomials in $r$. This allows us
   to define a deformation, or $q$-analogue, of these algebras which
   depends on a parameter~$q$. We determine the irreducible representations
   of all of these algebras and give a basis for their radicals. Finally,
   we show that the direct sum of cyclotomic Solomon algebras is
   canonically isomorphic to a concatenation Hopf algebra.
\end{abstract}
\maketitle

\section{Introduction}
In a seminal paper~\cite{Solomon:descent}, Solomon showed that the group
algebra of any finite Coxeter group has a remarkable subalgebra, the
\emph{Solomon descent algebra}. In this paper we construct a similar
subalgebra of the complex reflection group of type $G(r,1,n)$ and show that
this algebra shares many of the properties of the Solomon descent algebras.

Solomon showed that each descent algebra has a distinguished basis for
which he gave an explicit description of the structure constants. This
distinguished basis is given by the sums of the distinguished coset
representatives of the parabolic subgroups. Solomon gave a basis for the
radical of the descent algebra and he constructed a natural homomorphism
from the descent algebra into the parabolic Burnside ring of the associated
Coxeter group.  As a consequence, it follows that the irreducible
representations of the Solomon descent algebras are all one dimensional and
that, in characteristic zero, they are naturally indexed by the conjugacy
classes of the parabolic subgroups.

There has been an explosion of research into the descent algebras of
Coxeter groups
since Solomon discovered them; see, for example,~\cite{AtkPfiWill,BBHT,Bidigare:thesis,BHS:SolomonSym,BlesLaue,BonnafePfeiffer:descent,Pfeiffer:SolomonQuiver}.
The study of the Solomon descent algebras of the symmetric groups has been
even more intense because of connections between these algebras and free
Lie algebras, $0$-Hecke algebras, non-commutative and quasi-symmetric
functions~\cite{Atk:Descent,GarsiaReutenauer,GKLLTH,MalvenutoReutenauer:descent}, the
representation theory of the symmetric group, and card shuffling and
associated random walks~\cite{BayerDiaconis, Fulman}. 

The algebra that we construct in this paper is in many ways a natural
generalization of the Solomon algebra of the symmetric groups.  The
\emph{cyclotomic Solomon algebra} $\Sol(\G)$ is a subalgebra of the group
algebra of the complex reflection group~$\G$ of type $G(r,1,n)$.
Like Solomon, we define our algebra to be the subalgebra of the group
algebra of $\G$ with basis the `distinguished' coset representatives of a
natural class of subgroups of~$\G$. It turns out that many natural choices
of subgroups, and coset representatives for these subgroups, do not yield a
subalgebra of the group algebra (see Remark~\ref{not Mak}). We show,
however, that with respect to the `right' length function, the sums of the
minimal length coset representatives of the \textit{standard reflection
subgroups} of~$\G$ give rise to a subalgebra of $\Z\G$ which is free of
rank~$2\cdot3^{n-1}$. We give an explicit formula for the structure
constants for this basis which is similar to Solomon's formula for the
structure constants of the descent algebra of the symmetric group $\Sym_n$.

One surprising feature of the cyclotomic Solomon algebras $\Sol(\G)$ is
that the structure constants of these algebras for $n\ge0$ are polynomials
in~$r$ which are independent of~$n$. As a consequence, these algebras admit
a simultaneous deformation $\Sol_q(n)$ which depends on a parameter~$q$.
For fixed~$n\ge0$, we show that the algebras $\Sol_q(n)$ are free of rank
$2\cdot3^{n-1}$. We construct and classify the irreducible representations
of these algebras over an arbitrary field, and hence give a basis for the
radical of~$\Sol_q(n)$.

A remarkable result of Gessel~\cite{Gessel:pPartitions} shows that there is
a natural duality between the Hopf algebra of quasi-symmetric functions and
the descent Hopf algebra.  This led Malvenuto and
Reutenauer~\cite{MalvenutoReutenauer:descent} to show that the direct sum
of these algebras under the shuffle (or convolution) product is a Hopf
subalgebra of the Hopf algebra of permutations. This Hopf algebra is dual
to the Hopf algebra of quasi-symmetric functions and it is isomorphic to
the Hopf algebra of non-commutative symmetric functions \cite{GKLLTH}.
These results are important because they relate the coproduct of the
quasi-symmetric functions with the product in the descent algebras.  

Baumann and Hohlweg\cite{BaumannHohlweg:coloredhopf} showed that there is a
similar Hopf algebra structure under the shuffle product on the space
$\ColG=\bigoplus_{n\ge0}\Z\G$ of coloured permutations.  We prove that the
direct sum of the cyclotomic Solomon algebras
$\Sol(r)=\bigoplus_{n\ge0}\Z\Sol(\G)$ is  a Hopf subalgebra of~$\ColG$.
We show that $\Sol(r)$ is a concatenation Hopf algebra and that
$\Sol(r)$ has a second bialgebra structure which has the same coproduct as
$\ColG$ but where the product map is induced by group multiplication. We
expect that the Hopf algebra~$\Sol(r)$ is dual to the Hopf algebra of
quasisymmetric functions of type~$B$ considered by Hsiao and
Petersen~\cite{HsiaoPetersen:QuasiB}.

Different generalizations of the Solomon algebras have been considered by
other authors, the most striking of which are the Mantaci-Reutenauer
algebras~\cite{MantaciReutenauer:descent}. It is natural to ask whether the
cyclotomic Solomon algebras and the Mantaci-Reutenauer algebras are
isomorphic, at least for type $B_n$, since they are both free of
rank~$2\cdot3^{n-1}$. We show in Remark~\ref{not Mak} that, in general,
these two algebras are not isomorphic. Example~\ref{character} shows
that, in stark contrast to the Solomon descent algebra, there is no map
from $\Sol(\G)$ into the character ring of~$\G$.

This paper is organized as follows. In the second section we introduce the
complex reflection groups $\G$ and set our notation. In section~3 we define
and classify the standard reflection subgroups of~$\G$ and section~4 shows
that every coset of a reflection subgroup has a unique element of minimal
length. Sections~4 and~5 give combinatorial descriptions of the coset and
double representatives of the reflection subgroups. This combinatorics
turns out to be closely related to the structure constants of the
cyclotomic Solomon algebras, which are finally introduced in section~6.
The first main result of the paper, Theorem~\ref{cyc Solomon}, determines
the structure constants of the cyclotomic Solomon algebras, hence showing
that they are in fact subalgebras of~$\G$. In section~7 we investigate the
`generic' cyclotomic Solomon algebras and in section~8 we construct and
classify the irreducible representations of the cyclotomic algebras and
their deformations. In section~9 we show that the direct sum of the
cyclotomic algebras gives rise to a concatenation Hopf algebra which is a
Hopf subalgebra of the Hopf algebras of coloured permutations. Finally, in
section~10 we give a second combinatorial interpretation of the structure
constants of the cyclotomic Solomon algebras. We use this to show that the
direct sum of the cyclotomic Hopf algebras comes equipped with a second
bialgebra structure which has the same coproduct but where the product map
is induced by group multiplication.

\section{Complex reflection groups of type $G(r,1,n)$}

This paper is concerned with certain subalgebras of the group algebra
of the complex reflection groups of type $G(r,1,n)$, in the
Shephard--Todd classification of the finite subgroups of
$\operatorname{GL}_n(\C)$ which are generated by (pseudo) reflections. In this
section we introduce these groups and study a length function on
them.

Fix positive integers $r$ and $n$. The complex reflection group of
type $G(r,1,n)$ is the group $\G$ which is generated by
elements $s_0, s_1,\dots,s_{n-1}$ subject to the relations
\begin{align*}
      s_0^r &= 1 = s_i^2&
      s_0s_1s_0s_1 &= s_1s_0s_1s_0\\
      s_is_j&=s_js_i, &
      s_is_{i+1}s_i&=s_{i+1}s_is_{i+1},
\end{align*}
where $1\le i\le j-1\le n-1$. This presentation is very similar to the
presentation of a Coxeter group; indeed,  if $r\leq 2$ then $\G$ is a
Coxeter group. Accordingly, we encode this presentation in the following
``cyclotomic Dynkin diagram'':
\vskip 0in 
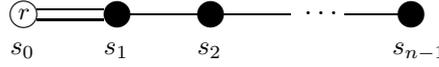
\begin{picture}(400,50)
\put(100,30){\circle{10}}
\put(98,28){\small{$r$}}
\put(105,32){\line(1,0){25}}
\put(105,28){\line(1,0){25}}
\put(135,30){\circle*{10}}
\put(140,30){\line(1,0){25}}
\put(170,30){\circle*{10}}
\put(175,30){\line(1,0){25}}
\put(205,30){\large{$\ldots$}}
\put(220,30){\line(1,0){25}}
\put(245,30){\circle*{10}}
\put(95,15){$s_0$}
\put(130,15){$s_1$}
\put(165,15){$s_2$}
\put(238,15){$s_{n-1}$}
\end{picture}

The node labeled by $r$ indicates that the generator $s_0$ has
order~$r$; otherwise, this graph gives the presentation of $\G$ in exactly
the same way as a Dynkin diagram gives the presentation of the
corresponding Coxeter group.  If $r=1$ then $G_{1,n}$ is isomorphic to the
symmetric group of degree~$n$. 

From the presentation of $\G$ it is evident that there is a
homomorphism from the symmetric group $\Sym_n$ into $\G$ which is
determined by mapping each transposition $(i,i+1)$ to $s_i$, for
$i=1,\dots,n-1$. In fact, this map is injective so we can -- and
do --  identify $\Sym_n$ with the subgroup $\<s_1,\dots,s_{n-1}\>$ via
this homomorphism. 

The symmetric group $\Sym_n$ acts on $\{1,2,\dots,n\}$ from the right.
We write this action exponentially. Thus,  $w\in\Sym_n$ sends the
integer $i$ to $i^w$, for $1\le i\le n$. 

Define $t_1=s_0$ and $t_{i+1}=s_it_is_i$, for $1\le i<n$.  Using the
relations it is easy to see that $t_it_j=t_jt_i$, for all $i,j$. It
follows that the subgroup $T=\<t_1,\dots,t_n\>$ is abelian and,
further, one can show that $T\cong(\Z/r\Z)^n$. It is easy to see
that
\begin{equation}\label{T left}
  t_i\, w = w\,t_{i^w},
 \quad\text{ for all $w\in\Sym_n$ and $1\le i\le n$,}
\end{equation}
Hence,~$T$ is a normal subgroup of $\G$. With a little more work we obtain
the following description of $\G$ as an (internal) semidirect product, or
wreath product:
\begin{equation}\label{wreath}
\G = T\rtimes\Sym_n = \<s_0\>\wr\<s_1,\dots,s_{n-1}\>
                    \cong (\Z/r\Z)\wr\Sym_n.
\end{equation}
Let 
$\Z_r^n=\set{\alpha=(\alpha_1,\dots,\alpha_n)\in\Z^n|0\le\alpha_i<r}$.
For $\alpha\in\Z_r^n$ let 
$t^\alpha=t_1^{\alpha_1}\dots t_n^{\alpha_n}$. Then, as a set,
$\G=\set{t^\alpha w|\alpha\in\Z_r^n\text{ and }w\in\Sym_n}$ and
$|\G|=r^nn!$.

Let $\Pi=\Pi_{r,n}=\{t_1,\dots,t_n,s_1,\dots,s_{n-1}\}$. Then
$\Pi$ generates $\G$ because $\{s_0=t_1,s_1,\dots,s_{n-1}\}$
generates~$\G$.
\begin{Defn}
The \textbf{$\Pi$--length function} on $\G$ is the 
function $\len=\len_\Pi\map{\G }\N$ given by
$\len(g)=\min\set{k\ge0|g=r_i\dots r_k, \text{ for some } r_i\in\Pi}$.
\end{Defn}

\begin{Remark}\label{length remark}
Let $S_0=\{s_0,s_1,\dots,s_{n-1}\}$. Bremke and
Malle~\cite{BM:redwds} have studied the length function
$\len_0\map{\G }\N$ which is defined by 
$$\len_0(g)=\min\set{k\ge0|g=r_i\dots r_k, \text{ for some } r_i\in S_0}.$$
By definition, $\len(g)\le\len_0(g)$, for all $g\in\G$. Furthermore,
it is not hard to see that $\len(g)\equiv\len_0(g)\pmod 2$. Moreover,
if $w\in\Sym_n$ then
$$\len(w)=\len_0(w)=\#\set{(i,j)|1\le i<j\le n\text{ and } i^w>j^w}.$$ 
(The last equality is well--known; see, for example,
\cite[Prop.~1.3]{M:ULect}.)
Hence, Proposition~\ref{T length} below gives an effective way of computing
the $\Pi$--length function on~$\G$.
\end{Remark}

For $\alpha=(\alpha_1,\dots,\alpha_n)\in\Z_r^n$ we set 
$|\alpha|=\alpha_1+\dots+\alpha_n$.

\begin{Prop}\label{T length}
Suppose that $\alpha\in\Z_r^n$ and $w\in\Sym_n$. Then
$\len(t^\alpha w)=|\alpha|+\len(w)$.
\end{Prop}

\begin{proof} By definition $\len(t^\alpha w)\le|\alpha|+\len(w)$.
    Conversely, suppose that $t^\alpha w=r_1\dots r_k$, for some
    $r_i\in\Pi$. Using (\ref{T left}) we can move each
    $t_i\in\{r_1,\dots,r_k\}$ to give a new word in which all of the
    elements of $T$ appear on the left. As every element of $\G$ can be
    written uniquely in the form $t^\beta v$, for
    $\beta\in\Z_r^n$ and $v\in\Sym_n$, this new word must be
    $t^\alpha w$.  By (\ref{T left}), this rewriting process does not
    increase the $\Pi$--length of the word, however, it may decrease the
    $\Pi$--length if some cancellation occurs. Hence,
    $k\ge|\alpha|+\len(w)$, completing the proof.
\end{proof}

\begin{Cor}\label{length} Suppose that $\alpha\in\Z_r^n$ and that
$w\in\Sym_n$. Then
\begin{align*}
\len(t_j\cdot t^\alpha w) &=\begin{cases}
      \len(t^\alpha w)+1,&\text{if }\alpha_j<r-1,\\
      \len(t^\alpha w)-r+1,&\text{if }\alpha_j=r-1,
\end{cases}
\end{align*}
for $1\le j\le n$ and $\len(s_i\cdot t^\alpha w) = |\alpha|+\len(s_iw)$,
for $1\le i<n-1$.
\end{Cor}

Note that $t^\alpha w\cdot t_j=t_{j^{w^{-1}}}\cdot t^\alpha w$ by (\ref{T
left}) and $\len(t^\alpha w\cdot s_i)=|\alpha|+\len(ws_i)$, for $1\le i<n-1$
and $1\le j\le n$. Hence, Corollary~\ref{length} can be used to compute
$\len(g\cdot t^\alpha w)$ and $\len(t^\alpha w\cdot g)$, for any $g\in\G$.

It is sometimes convenient to describe $\G$ combinatorially as a set of `words'.
Fix a primitive $r^{\text{th}}$ root of unity $\zeta=\exp(2\pi i/r)\in\C$ and set
$$\n=\{1,2,\dots,n\}\qquad\text{and}\qquad
  \rn=\set{m\zeta^i| m\in\n\text{ and }0\le i<r}.$$ 
Recall that if $z\in\C$ then $|z|$ is the complex modulus of $z$. In particular, if $m\zeta^i\in\rn$ then $|m\zeta^i|=m$.  Define a 
\textbf{word} in $\rn$ to be an element of the set
$$\Words=\set{\underline\omega=(\omega_1,\dots,\omega_n)|
\omega_i\in\rn\text{ and }\{|\omega_1|,\dots,|\omega_n|\}=\n}.$$ 
If $\omega=(\omega_1,\dots,\omega_n)$ is a word then we abuse notation and
write $\omega=\omega_1\dots\omega_n$.

There is a faithful right action of~$\G$ on~$\Words$ given by
$$\omega_1\dots\omega_n\cdot t^\alpha w
=\zeta^{\alpha_{1}}\omega_{1^w}\dots\zeta^{\alpha_{n}}\omega_{n^w},$$
for $\alpha\in\Z_r^n$ and $w\in\Sym_n$. Consequently, there is
a natural bijection $\G\to\Words$ given by 
$t^\alpha w\mapsto1\dots n\cdot t^\alpha w$, so that $|\Words|=r^nn!=|\G|$.
Thus, we have described the regular representation of~$\G$ as the
permutation representation on the set of words $\Words$. Equivalently, $\G$
is the group of permutations of $\rn$ such that $(m\zeta^i)^g=m^g\zeta^i$,
for all $m\in\n$, $0\le i<r$ and $g\in\G$.

\section{Reflection subgroups}

Recall that $\Pi=\{t_1,\dots,t_n,s_1,\dots,s_{n-1}\}$. In this section
we define the reflection subgroups of $\G$ and show that every coset
of a reflection subgroup contains a unique element of minimal
$\Pi$--length.

\begin{Defn} A \textbf{$($standard$)$ reflection subgroup} of $\G$ is a
  subgroup which is generated by a subset of $\Pi$.
\end{Defn}

Geometrically, a reflection subgroup of $\G$ should be any subgroup which
is generated by elements which act by (pseudo) reflections in the
reflection representation of $\G$. All of the elements of $\Pi$ act as
reflections in the reflection representation of~$\G$, so every standard
reflection subgroup is a reflection subgroup in this geometric sense. If
$r>2$ then it is not difficult to see that there are `geometric reflection
subgroups' of $\G$ which are not standard reflection subgroups.

If $J\subseteq\Pi$ let $G_J=\<J\>$ be the corresponding (standard)
reflection subgroup of~$\G$. This notation is inherently
ambiguous because it can happen that $G_J=G_K$ even though $J\ne K$, for
$J,K\subseteq\Pi$. For example, $G_\Pi=\G=G_{S_0}$ (recall that
$S_0=\{s_0,s_1,\dots,s_{n-1}\}$), and yet $\Pi\ne S_0$ if $n>1$. We start our
study of the reflection subgroups by resolving this ambiguity.

A \textbf{composition} of $n$ is a sequence $\mu=(\mu_1,\dots,\mu_k)$ of
positive integers which sum to~$n$. A \textbf{signed composition} of $n$ is
a sequence of non--zero integers $\mu=(\mu_1,\dots,\mu_k)$ such that
$|\mu|=|\mu_1|+\dots+|\mu_k|=n$. Let $\SComp$ be the set of signed compositions
of~$n$ and let $\Comp$ be the set of compositions of $n$. Then
$\Comp\subseteq\SComp$.

If $\mu=(\mu_1,\dots,\mu_k)\in\SComp$ let
$\mu^+=(|\mu_1|,\dots,|\mu_k|)$ and $-\mu=(-\mu_1,\dots,-\mu_k)$. Then
$\mu^+\in\Comp$ is a composition of~$n$ and $-\mu\in\SComp$. We set
$|\mu|^+=\frac12\sum_{i=1}^k(\mu_i^++\mu_i)$, so that~$|\mu|^+$ is the sum of
the positive parts of $\mu$. Similarly, let
$|\mu|^-=\frac12\sum_{i=1}^k(\mu_i^+-\mu_i)$ be the absolute value of the
sum of the negative parts of $\mu$. Then $|\mu|=|\mu|^-+|\mu|^+=n$.
Finally, set $\bar\mu_0=0$ and $\bar\mu_i=|\mu_1|+\dots+|\mu_i|$, for $i\ge
1$.  

\begin{Defn}\label{J mu}
  Suppose that $\mu=(\mu_1,\dots,\mu_k)\in\SComp$ is a signed composition. 
  Define
  $$\Pi_\mu=\bigcup_{1\le i\le k}\{s_{\bar\mu_{i-1}+1},\dots,s_{\bar\mu_i-1}\}
  \cup\bigcup_{\substack{1\le i\le k\\\mu_i>0}}
  \{t_{\bar\mu_{i-1}+1},\dots,t_{\bar\mu_i}\}.$$
  Then $\Pi_\mu\subseteq\Pi$ so we set $G_\mu=G_{\Pi_\mu}$.
\end{Defn}

Let $S= \{ s_1, \ldots, s_{n-1}\}\subseteq\Pi$. Suppose that $\mu\in\SComp$. Then
$\Pi_\mu\subseteq S$ if and only if $-\mu\in\Comp$.
In general, $\Pi_\mu\subseteq\Pi$ and the reflection subgroup 
$G_\mu$ is conjugate to the reflection subgroup
$$\prod_{\mu_i>0}G_{r,\mu_i}\,\times\,
\prod_{\mu_j<0}\Sym_{-\mu_j}$$
of $\G$. Moreover, $\set{G_\mu|\mu\in\SComp}$ is the complete set of
reflection subgroups of~$\G$.

\begin{Prop}\label{reflection subgroups}
  Suppose that $n\ge 1$, $r\geq 2$ and that $J\subseteq\Pi$. Then 
  $G_J=G_\mu$, for a unique signed composition $\mu\in\SComp$. Consequently,
  $\G$ has $2\cdot3^{n-1}$ distinct reflection subgroups.
\end{Prop}

\begin{proof} 
  We prove both statements in the Proposition by induction on~$n$. If
  $n=1$ then $G_\emptyset=G_{(1)}$ and $G_\Pi=G_{(-1)}$ are the only
  reflection subgroups of $G_{r,1}$ so the Proposition holds. In
  particular, $G_{r,1}$ has $|\SComp[1]|=2$ reflection subgroups.

  Suppose then that $n>1$ and observe that
  $\Pi_{r,n}=\Pi_{r,n-1}\cup\{s_{n-1},t_n\}$. Let $G'=G_{r,n-1}$, which we
  consider as a subgroup of $\G$ in the natural way. By induction on $n$
  every reflection subgroup of $G'$ is of the form
  $G'_\mu=(G_{r,n-1})_\mu$, for some
  $\mu\in\SComp[n-1]$. 
  
  Fix $J\subseteq\Pi$. Then $G_J\cap G'$ is a reflection subgroup of $G'$,
  so that $G_J\cap G'=G_\mu$, for some $\mu\in\SComp[n-1]$. Now, 
  $t_{n-1}\in G'_\mu$ if and only if $\mu_k>0$, so one can check that
  $$\<G_{(\mu_1,\dots,\mu_k)}',s_{n-1},t_n\>=\begin{cases}
    \<G_{(\mu_1,\dots,\mu_{k-1},\mu_k)}',s_{n-1}\>,&\text{if }\mu_k>0,\\
    \<G_{(\mu_1,\dots,\mu_{k-1},-\mu_k)}',s_{n-1}\>,&\text{if }\mu_k<0.
  \end{cases}$$
  Consequently, $G_J$ is equal to either $G'_\mu$, $\<G'_\mu,t_n\>$ or 
  $\<G'_\mu,s_{n-1}\>$. Therefore,
  $$G_J=\begin{cases}
    G'_{(\mu_1,\dots,\mu_k)}=G_{(\mu_1,\dots,\mu_k,-1)},
       &\text{if } s_{n-1},t_n\notin G_J\\
    \<G'_{(\mu_1,\dots,\mu_k)},t_n\>=G_{(\mu_1,\dots,\mu_k,1)},
      &\text{if } s_{n-1}\notin G_J\text{ and }t_n\in G_J\\
    \<G'_{(\mu_1,\dots,\mu_k)},s_{n-1}\>=G_{(\mu_1,\dots,\mu_k+\epsilon_k)},
       &\text{if } s_{n-1},t_n\in G_J\\
  \end{cases}$$
  where $\epsilon_k=1$ if $\mu_k>0$ and $\epsilon_k=-1$ if $\mu_k<0$.
  Hence, the reflection subgroups of $G_{r,n}$ are naturally indexed by
  the signed compositions of $n$. Moreover, by
  (\ref{wreath}) the subgroups of~$\G$ arising this way for different
  $\nu\in\SComp[n-1]$ are all distinct. Consequently, by induction, 
  $G_{r,n}$ has $3|\SComp[n-1]|=2\cdot 3^{n-1}$ reflection subgroups.
\end{proof}

It follows from the definitions and Proposition~\ref{reflection subgroups}
that $\Pi_\mu$ is the unique maximal subset of $\Pi$ (under inclusion)
which generates the reflection subgroup $G_\mu$. In contrast, if
$\mu\in\SComp$ then the reader can check that there are
$\prod_{i:\mu_i>0}\mu_i$ distinct minimal subsets of~$\Pi$ which generate
$G_\mu$. Thus, the (minimal) subsets of~$\Pi$ which generate the reflection
subgroups are, in general, not unique.

\section{Distinguished coset representatives.}
In this section we describe, both algebraically and combinatorially, a set
of `distinguished' coset representatives for the reflection
subgroups of~$\G$.

Fix a composition $\lambda=(\lambda_1,\dots,\lambda_k)$ of $n$. Then
$\Sym_\lambda=\Sym_{\lambda_1}\times\dots\times\Sym_{\lambda_k}$ is a
parabolic, or Young subgroup of $\Sym_n$. According to our conventions
$\Sym_\lambda=G_{-\lambda}$, so $\Sym_\lambda$ is a reflection subgroup
of~$\G$. Let 
$$\D_\lambda=\set{d\in\Sym_n| \len(d)\le\len(w)\text{ for all }w\in\Sym_\lambda d}.
$$
Then, as is well--known, $\D_\lambda$ is a complete set of right coset
representatives for $\Sym_\lambda$ in~$\Sym_n$. Moreover, if
$d\in\D_\lambda$ then $d$ is the unique element of minimal length in the
coset $\Sym_\lambda d$; see, for example,
\cite[Prop.~2.1.1]{GeckPfeiffer:book}. It is
not hard to see that $T\D_\lambda$ is a complete set of minimal length
coset representatives for $G_{-\lambda}=\Sym_\lambda$ in~$\G$. We want to
generalize this observation to all reflection subgroups.

Recall that if $\mu=(\mu_1,\dots,\mu_k)\in\SComp$ then
$\mu^+=(|\mu_1|,\dots,|\mu_k|)$ is a composition of~$n$. Consulting the
definitions, $G_\mu\cap\Sym_n=\Sym_{\mu^+}$. Similarly, define 
\begin{align*}
  T_\mu&=G_\mu\cap T 
  =\<t_i\mid t_i\in G_\mu\>\\
  &=\<t_i\mid \bar\mu_{j-1}<i\le\bar\mu_j\text{ for some $j$ with }\mu_j>0\>.
\end{align*}
Then, $T_\mu\cong(\Z/r\Z)^{|\mu|^+}$.

With this notation, (\ref{wreath}) gives the following description of
$G_\mu$ as a semidirect product of $T_\mu$ and $\Sym_{\mu^+}$.

\begin{Lemma}\label{semidirect} 
Suppose that $\mu\in\SComp$. Then $G_\mu=T_\mu\rtimes\Sym_{\mu^+}$.
\end{Lemma}

Since $T\cong(\Z/r\Z)^n$ is an abelian group, every subgroup of~$T$
is a normal subgroup of~$T$. In particular, if $G_\mu$ is a
reflection subgroup of $\G$ then $T_\mu$ is normal in~$T$ and 
$T/T_\mu\cong (\Z/r\Z)^{|\mu|^-}\cong T_{-\mu}$. Further,
$T_\mu T_{-\mu}=T=T_{-\mu}T_\mu$, for all $\mu\in\SComp$.

Mimicking the definition of $\D_\mup$ we have:

\begin{Defn}Suppose that $\mu\in\SComp$. Set
  $$\E_\mu = \set{e\in\G|\len(e)\le\len(g)\text{ for all }g\in G_\mu e}.$$
\end{Defn}

We can now prove the main result of this section which shows that
$\E_\mu$ is a (distinguished) set of coset representatives for
$G_\mu$ in~$\G$.

\begin{Thm}\label{Dmu} 
  Suppose that $\mu\in\SComp$. Then $\E_\mu=T_{-\mu}\times\D_\mup$ and
  $\E_\mu$ is a complete set of right coset representatives for $G_\mu$
  in~$\G$.
\end{Thm}

\begin{proof} We first show that $T_{-\mu}\times\D_\mup$ is a complete set
  of coset representatives for~$G_\mu$ in~$\G$. Suppose that
  $t^\alpha w\in\G$, where $\alpha\in\Z_r^n$ and $w\in\Sym_n$.
  Define  $\beta=(\beta_1,\dots,\beta_n)\in\Z_r^n$ by
  $$\beta_i=\begin{cases}
    \alpha_i,&\text{if }t_i\notin G_\mu\iff t_i\in T_{-\mu},\\
    0,&\text{if }t_i\in G_\mu\iff t_i\notin T_{-\mu}.
  \end{cases}$$
Then, by definition, $t^\beta\in T_{-\mu}$. Moreover,
$G_\mu t^\alpha w= G_\mu t^\beta w$ and $\len(t^\alpha
w)\ge\len(t^\beta w)$, with equality if and only if $\alpha=\beta$.

Write $w=vd$, where $v\in\Sym_{\mup}$ and $d\in\D_\mup$. Let $\gamma=\beta
v=(\beta_{1^v},\dots,\beta_{n^v})$. Then $t^\beta v=vt^\gamma$, by (\ref{T
left}), so that $t^\gamma=v^{-1}t^\beta v\in T_{-\mu}$ since $\Sym_\mu$
centralizes $T_{-\mu}$. Consequently, $G_\mu t^\alpha w=G_\mu t^\gamma d$,
where $t^\gamma\in T_{-\mu}$ and $d\in\D_\mup$. However, by
Lemma~\ref{semidirect}, $$[\G
:G_\mu]=[T:T_{-\mu}]\cdot[\Sym_n:\Sym_{\mu^+}]
=\#(T_{-\mu}\times\D_\mup).$$ Therefore, $T_{-\mu}\times\D_\mup$ is a
complete set of right coset representatives for $G_\mu$ in $\G$.

It remains to prove that $\E_\mu=T_{-\mu}\times\D_\mup$. Suppose that, as
above, we have $G_\mu t^\alpha w=G_\mu t^\gamma d$, for $\alpha\in\Z_r^n$,
$w\in\Sym_n$, $t^\gamma\in T_{-\mu}$ and $d\in\D_\mup$. The argument of the
first paragraph shows that $\len(t^\gamma d)\le\len(t^\alpha w)$ with
equality if and only if $t^\alpha\in T_{-\mu}$ and $w\in\D_\mup$. That is,
if and only if $\alpha=\gamma$ and $w=d$.  Hence,
$\E_\mu=T_{-\mu}\times\D_\mup$ as claimed.
\end{proof}

Theorem~\ref{Dmu} shows that every coset of a reflection subgroup
contains a unique element of minimal $\Pi$--length. We call $\E_\mu$ the
set of \textbf{distinguished coset representatives} for $G_\mu$
in~$\G$.

\begin{Example}\label{refl example} Suppose that $r\ge2$ and consider 
  $G_{r,2}=(\Z/r\Z)\wr\Sym_2$. Then
  $\Pi=\{t_1,t_2,s_1\}$ and $G_{r,2}$ has six reflection subgroups. The
  following table describes these groups and the corresponding sets of
  distinguished right coset representatives.
  $$\begin{array}{*4c}
    \mu   & G_\mu  & \Pi_\mu & \E_\mu \\[2pt]\toprule
    (-1,-1) & 1      &  \emptyset        & T\times\Sym_2\\
    (1,-1)& \set{t_1^k|0\le k<r} &\{t_1\}&\<t_2\>\times\Sym_2\\
    (-1,1)& \set{t_2^k|0\le k<r} &\{t_2\}&\<t_1\>\times\Sym_2\\
    (1,1)& T      & \{t_1,t_2\}         & \Sym_2\\
    (-2)    & \Sym_2 &   \{s_1\}           & T\\
    (2)& T\rtimes\Sym_2 & \{t_1,t_2,s_1\}   & 1
  \end{array}$$
  For each reflection subgroup we have given the factorization of $\E_\mu$
  from Theorem~\ref{Dmu}. Observe that the
  reflection subgroups do not depend in a crucial way on~$r$ and that
  $|\E_\mu|=|\G|/|G_\mu|$ is a polynomial in~$r$, for $\mu\in\SComp$
  (and $r\ge2$).
\end{Example}

We now give combinatorial interpretations of the set of distinguished
coset representatives $\E_\mu$, for $\mu\in\SComp$, which is
similar to the description of $\D_{\mu^+}$ in terms of row standard
tableaux (see \cite[Prop.~3.3]{M:ULect}).

Fix a composition $\lambda\in\Comp$. The \textbf{diagram} of $\lambda$ is the set
$$[\lambda]=\set{(i,j)\in\N^2|1\le j\le\lambda_i\text{ and }1\le i\le \ell(\lambda)}.$$
Here $\ell(\lambda)$ is the number of non--zero parts of $\lambda$.
We think of $[\lambda]$ as being an array of boxes in the plane.

Now suppose that $\mu\in\SComp$. A \textbf{$\mu$--tableau} is a map
$\t\map{[\mu^+]}\rn$. We identify a $\mu$--tableau with a diagram for
$\mup$ which is labeled by elements of~$\rn$. If $\t$ is a $\mu$--tableau
let $|\t|$ be the tableau obtained by taking the complex modulus of the
entries in $\t$; that is, $|\t|(x)=|\t(x)|$, for all $x\in[\mu^+]$.

\begin{Example}\label{tableaux}
Let $\mu=(2,-3,1,-1)$. Then four $\mu$--tableaux are:
$$\Tab(1&2\cr3&4&5\cr6\cr7)\qquad
  \Tab(1&2\cr3\zeta&4\zeta^2&5\zeta^3\cr6\cr7\zeta^4)\qquad
  \Tab(3&6\cr2\zeta&5\zeta^2&7\zeta^3\cr1\cr4\zeta)\qquad\text{and}\qquad
  \Tab(7&6\zeta\cr2\zeta&5\zeta^2&3\zeta^3\cr1\cr4\zeta).
$$
\end{Example}

As remarked at the end of section~2 we can think of $\G$ as the group of
permutations of~$\rn$ such that $(m\zeta^i)^g=m^g\zeta^i$, for all
$m\zeta^i\in\rn$ and all $g\in\G$. Consequently, $\G$ acts on the set of
$\mu$--tableaux by composition of maps. Thus, if~$\t$ is a $\mu$--tableau
and $g\in\G$ then $\t^g$ is the tableau with $\t^g(x)=\t(x)^g$, for
$x\in[\mu^+]$. 

Let $\t^\mu$ be the $\mu$--tableau which has the numbers $1,\dots,n$
entered in order, from left to right and then top to bottom, along the
rows of $[\mu^+]$. The first $\mu$--tableau in
Example~\ref{tableaux} is~$\t^\mu$ when $\mu=(2,-3,1,-1)$.

So far none of the combinatorial definitions above distinguish between
compositions and signed compositions. We now single out a set of
$\mu$--tableaux that are in bijection with~$\E_\mu$. First, define a total
order $\preceq$ on $\rn$ by declaring that
$a\zeta^i\preceq b\zeta^j$ if $a<b$, or $a=b$ and $i>j$. Then
$$\zeta^{m-1}\preceq\zeta^{m-2}\preceq\dots\preceq\zeta\preceq1
   \preceq2\zeta^{m-1}\preceq\cdots\preceq2\preceq\dots
   \preceq n\zeta^{m-1}\preceq \cdots \preceq n.
$$

\begin{Defn}\label{row standard tableaux}
Suppose that $\mu\in\SComp$. A $\mu$--tableau $\t$ is 
\textbf{row standard} if it satisfies the following three conditions:
\begin{enumerate}
\item The set of entries in the tableau $|\t|$ is $\{1,\dots,n\}$.
\item The entries in row $i$ of $\t$ belong to $\{1,\dots,n\}$ 
  whenever $\mu_i>0$.
\item In each row the entries of $\t$ appear, from left to right, in 
    increasing order with respect to $\preceq$.
\end{enumerate}
\end{Defn}

For example, the first three of the $(2,-3,1,-1)$--tableaux in
Example~\ref{tableaux} are row standard. The last tableau in this
example is not row standard because it fails conditions (b) and (c).

The action of $\G$ on the set of $\mu$--tableau which satisfy condition~(a)
of Definition~\ref{row standard tableaux} gives a realization of the regular
representation of~$\G$. Consequently, the map $g\mapsto\t^\mu g$, for
$g\in\G$, is a bijection from $\G$ to the set of these $\mu$--tableaux.  If
$\t$ is such a $\mu$--tableau let $d_\t$ be the unique element of $\G$ such
that $\t=\t^\mu d_\t$.

\begin{Prop}\label{row standard}
Suppose that $\mu\in\SComp$. Then
$$\E_\mu=\set{d_\t|\t\text{ is a row standard $\mu$--tableau}}.$$
\end{Prop}

\begin{proof}
By definition, the orbit $\t^\mu G_\mu=\set{\t^\mu g|g\in G_\mu}$ of
$\t^\mu$ under $G_\mu$ consists of all those tableaux which can be
obtained by permuting the entries of each row of~$\t^\mu$ and
multiplying the entries in row~$i$ by a power of~$\zeta$ when
$\mu_i>0$. Consequently,~$\t^\mu$ is the unique row standard
$\mu$--tableaux in $\t^\mu G_\mu$, so that each right coset
of $G_\mu$ in~$\G$ contains a unique element $e$ such that $\t^\mu e$
is row standard. Now, $\E_\mu=T_{-\mu}\D_\mup$ by
Theorem~\ref{Dmu} and $T_{-\mu}$ acts on $\t^\mu$ by multiplying the
entries in row~$i$ by different powers of~$\zeta$ when $\mu_i<0$. If
$d\in\D_\mup$ then it is well--known that the entries in $\t^\mu d$
increase from left to right along each row; see, for example,
\cite[Prop.~3.3]{M:ULect}. Hence, $\t^\mu e$ is row standard whenever 
$e\in\E_\mu$. This completes the proof.
\end{proof}

In the case of the symmetric groups the set of distinguished coset
representatives can be described combinatorially in terms of `descents'.
Explicitly, if $w\in\Sym_n$ then its \textbf{descent set} is 
$$\Des(w)=\set{s\in S|\ell(sw)<\ell(w)}
         =\set{s_i|1\le i<n\text{ and }i^w>(i+1)^w}.$$
If $\mu$ is a composition of $n$ then the connection between
distinguished coset representatives and descents is that
\begin{equation}\label{descent}
  \D_\mu=\set{d\in\Sym_n|\Des(d)\subseteq S-\Pi_{-\mu}}.
\end{equation}

There is an analogous description of $\E_\mu$, for
$\mu\in\SComp$. If $\alpha\in\Z^n_r$ define the \textbf{colour} of~$t^\alpha$
to be the set $\Col(t^\alpha)=\set{t_i\in T|\alpha_i>0}.$
Then using Theorem~\ref{Dmu} it is easy to see that if
$\mu\in\SComp$ then
$$\E_\mu=\set{t^\alpha w\in\G|\Col(\alpha)\cup\Des(w)\subseteq\Pi-\Pi_{-\mu}}.$$
We remark that it is easy to rephrase this last statement
combinatorially in terms of words in~$\Words$.

\begin{Remark} It is easy to check that
$\E_\mu^{-1}=\D_\mu^{-1}\times T_{-\mu}$ is a
complete set of left coset representatives for $G_\mu$ in~$\G$.
Moreover, $e\in\E_\mu^{-1}$ if and only if $\len(e)\le\len(g)$ for all
$g\in eG_\mu$, so every left coset of $G_\mu$ in $\G$ contains a
unique element of minimal $\Pi$--length. 
\end{Remark}

\begin{Remark}\label{Mak}
Mak~\cite{Mak:parabolic} has shown that every coset of
a reflection subgroup contains a unique element of minimal length with
respect to the length function~$\len_0$ defined in
Remark~\ref{length remark}. Mak's set of coset representatives is
different from $\E_\mu$. Nonetheless, it does admit a factorization
which is similar to the factorization of~$\E_\mu$ given in
Theorem~\ref{Dmu}. To describe this if $\mu=(\mu_1,\dots,\mu_k)\in\SComp$ 
then set
$$\E_\mu'=\prod_{\substack{k\ge j\ge 1\\\mu_j<0}}
  \prod_{\bar\mu_j\ge i>\bar\mu_{j-1}}\Big[\{1\}\cup
  \set{s_{\bar\mu_{j-1}+1}s_{\bar\mu_{j-1}+2}\dots s_{i-1} t_i^k|1\le k<r}
  \Big]{\times}\D_\mu,$$
where the product is taken in order from left to right in terms of
decreasing values of~$i$. One can show that $\E_\mu'$ is Mak's set
of coset representatives for $G_\mu$ in~$\G$. As we will never
need this result we leave the proof to the reader.
\end{Remark}

\section{Double Coset representatives}
Our next aim is to describe the double cosets of reflection subgroups. 
In order to do this we first recall some well--known facts about the
symmetric group $\Sym_n$.  Suppose that~$\mu$ and~$\nu$ are compositions
of~$n$. Then $\Sym_\mu$ and $\Sym_\nu$ are Young, or parabolic, subgroups
of $\Sym_n$.  Set $\D_{\mu\nu}=\D_\mu\cap\D_\nu^{-1}$. Then $\D_{\mu\nu}$
is a complete set of $(\Sym_\mu,\Sym_\nu)$--double coset representatives
in~$\Sym_n$; see, for example, \cite[Prop.~4.4]{M:ULect}. Moreover, if
$d\in\D_{\mu\nu}$ then $d^{-1}\Sym_\mu d\cap\Sym_\nu$ is  a Young subgroup
of~$\Sym_n$; see, for example, \cite[Lemma~4.3]{M:ULect}.  Define $\mu
d\cap \nu$ to be the unique composition of~$n$ such that $\Sym_{\mu
d\cap\nu}=d^{-1}\Sym_\mu d\cap\Sym_\nu$. We remark that the composition
$\mu d\cap \nu$ can be determined by comparing the row stabilizers of the
tableaux $\t^\mu d$ and $\t^\nu$.

\begin{Lemma}\label{mud cap nu}
  Suppose that $\mu,\nu\in\SComp$ and $d\in\D_{\mup\nup}$. 
  Then $d^{-1}G_\mu d\cap G_\nu$ is a reflection subgroup of~$\G$.
\end{Lemma}

\begin{proof}The group $G_\nu$ consists of those elements of~$\G$ which act
  on~$\t^\nu$ by first multiplying each entry of row~$l$ by possibly
  different powers of $\zeta$, if $\nu_l>0$, and then permuting the entries
  in each row of the resulting tableaux. Similarly, the group $d^{-1}G_\mu
  d$ consists of those elements of~$\G$ which act on the row standard
  tableau~$\t^\mu d$ by multiplying each entry of row~$k$ by different
  powers of $\zeta$, if $\mu_k>0$, and then permuting the entries in each
  row. Consequently, the subgroup $d^{-1}G_\mu d\cap G_\nu$ is generated by
  the elements $\{s_i,t_j\}$, where~$i$ runs over those integers for which
  $i$ and $i+1$ are in the same row of $\t^\nu$ and in the same row
  of~$\t^\mu d$, and $j\in\n$ is in row~$l$ of $\t^\nu$ with~$\nu_l>0$ and
  $j$ is in row~$k$ of $\t^\mu d$ with~$\mu_k>0$ (cf.~the proof of
  \cite[Lemma~4.3]{M:ULect}).  Therefore,~$d^{-1}G_\mu d\cap G_\nu=G_\sigma$,
  where $\sigma$ is the unique signed composition such that
  $\sigmap=\mup d\cap\nup$ and $\sigma_i>0$ if and only if $\nu_j>0$ and
  $\mu_k>0$, where $\bar\sigma_i$ appears in row~$j$ of $\t^\nu$ and
  row~$k$ of~$\t^\mu d$.  \end{proof}

Suppose that $d\in\D_{\mup\nup}$, for $\mu,\nu\in\SComp$. Then
$d^{-1}\in\D_{\nup\mup}$, since $\D_{\nup\mup}=\D_{\mup\nup}^{-1}$.
Therefore, $G_\mu\cap dG_\nu d^{-1}$ is
also a reflection subgroup of $\G$.

\begin{Defn}\label{intersect scomp}
  Suppose that $\mu,\nu\in\SComp$ and $d\in\D_{\mup\nup}$. Then $\mu
  d\cap\nu=\nu\cap\mu d$ is the signed composition of~$n$ such that $G_{\mu
  d\cap\nu}=d^{-1}G_\mu d\cap G_\nu$ and $\mu\cap d\nu=d\nu\cap\mu$ is the
  signed composition such that $G_{\mu\cap d\nu}=G_\mu\cap dG_\nu d^{-1}$.
\end{Defn}

Note that the proof of Lemma~\ref{mud cap nu} gives a recipe for computing
$\mu d\cap\nu$. Note also that $\mu d\cap\nu=d^{-1}\mu\cap\nu$, for
$d\in\D_{\mup\nup}$ and $\mu,\nu\in\SComp$.

We now describe a set of $(G_\mu,G_\nu)$--double coset representatives. We
do this by generalizing the description of the double cosets of the Young
subgroups of the symmetric group in terms of row semistandard tableaux. 

\begin{Defn}Suppose that $\mu\in\SComp$. A $\mu$--tableau 
$\T\map{[\mup]}\rn$ is \textbf{row semistandard} if
\begin{enumerate}
  \item The entries in row $i$ of $\T$ belong to $\{1,\dots,n\}$ 
    whenever $\mu_i>0$.
  \item The entries in each row of $\T$ appear in weakly increasing order,
    from left to right, with respect to~$\preceq$.
\end{enumerate}
\end{Defn}

There is a map from the set of row semistandard tableau to the set of row
standard tableaux. To define this first observe that a row
semistandard $\mu$--tableau $\T$ determines a unique total order $<_\T$ on
$[\mup]$ where $x<_\T x'$, for $x,x'\in[\mup]$, if 
\begin{enumerate}
\item $|\T(x)|<|\T(x')|$, or 
\item $|\T(x)|=|\T(x')|$ and $x$ is in an earlier row of $[\mup]$ than~$x'$, or 
\item $|\T(x)|=|\T(x')|$ and $x$ and $x'$ are in the same row and~$x$ is
to the left of~$x'$. 
\end{enumerate}
Let $x_1<_\T\dots<_\T x_n$ be the nodes in
$[\mup]$. Then the $\mu$--tableau~$\T^*$ is defined by the requirements
that $|\T^*(x_i)|=i$ and $\arg\T^*(x_i)=\arg\T(x_i)$, for $1\le i\le n$.
(If $z\in\C$ is a complex number let $\arg z\in[0,2\pi)$ be its argument so
that $z=|z|\exp(i\arg z)$.) By construction, $\T^*$ is a row standard
$\mu$--tableaux. Moreover, it is easy to see that the map $\T\mapsto\T^*$
is injective. 

\begin{Defn}
Suppose that $\mu,\nu\in\SComp$ and let $\T$ be a $\mu$--tableau. Then 
$\T$ has \textbf{type $\nu$} if
    \begin{enumerate}
	    \item $|\nu_j|=\#\set{x\in[\mup]| |\T(x)|=j}$, for $j\ge1$.
      \item If $\nu_j>0$ then $\nu_j=\#\set{x\in[\mup]|\T(x)=j}$.
    \end{enumerate}
Let $\srssTab(\mu,\nu)=\set{\T\map{[\mup]}\rn|\text{$\T$ is row semistandard
 $\mu$--tableau of type $\nu$}}$.
If~$\mu$ and~$\nu$ are \emph{compositions} let
$\rssTab(\mu,\nu)=\set{\T\map{[\mup]}\n|\text{$\T$ is row semistandard
 $\mu$--tableau of type $\nu$}}$.
\end{Defn}

See Example~\ref{row semistandard} below for these definitions in action.

\begin{Point}[\protect{\cite[Prop.~4.4]{M:ULect}}]\label{double Sym}
  Suppose that $\mu,\nu\in\Comp$. Then 
  $$\D_{\mu\nu}=\set{d_{\T^*}|\T\in\rssTab(\mu,\nu)}$$
  is a complete set of $(\Sym_\mu,\Sym_\nu)$ double coset representatives
  in $\Sym_n$. Moreover, if $d\in\D_{\mu\nu}$ then $\ell(d)\le\ell(w)$, for
  all $w\in\Sym_\mu d\Sym_\nu$, with equality if and only if $w=d$.
\end{Point}

If $\t$ is a row standard tableau let $\nu(\t)'$ be the tableau
obtained by replacing each entry~$m\zeta^a$ in $\t$ with $k\zeta^{a'}$ if
$m$ appears in row~$k$ of $\t^\nu$, where $a'=0$ if $\nu_k>0$ and $a'=a$
otherwise. Now define $\nu(\t)$ to be the row semistandard tableau obtained
by reordering the entries in each row of~$\nu(\t)'$ so that they are in
increasing order. Then $\nu(\t)$ is a row semistandard tableau of
type~$\nu$.

For example, let $\nu=(2,-2,1)$ and 
$\t = \Tab(3&4\zeta^a &5\zeta^b\cr 1&2)$~, where $0\le a,b<r$. Then, by
definition, $\nu(\t)' =\Tab(2&2\zeta^a &3\cr 1&1)$\quad and \
$\nu(\t) =\Tab(2\zeta^a&2 &3\cr 1&1)$~.

\begin{Prop}\label{double}
Suppose that $\mu$ and $\nu$ are signed compositions of $n$ and let
$$\E_{\mu\nu}=\set{d_{\T^*}|\T\in\srssTab(\mu,\nu)}.$$
Then $\E_{\mu\nu}$ is a complete set of $(G_\mu,G_\nu)$~double coset
representatives in~$\G$.  Moreover, if~$e\in\E_{\mu\nu}$ then 
$\ell(e)\le\ell(g)$, for all $g\in G_\mu eG_\nu$.
\end{Prop}

\begin{proof}
By Proposition~\ref{row standard} the right cosets of $G_\mu$ in $\G$ are
naturally indexed by the row standard $\mu$--tableaux. Hence, the
$(G_\mu,G_\nu)$--double cosets are indexed by the $G_\nu$--orbits
of the row standard $\mu$--tableaux. Using the definitions it is easy to
see that two $\mu$--tableaux~$\s$ and $\t$ belong to the same
$G_\nu$--orbit if and only if $\nu(\s)=\nu(\t)$. Moreover, if $\t$ is row
standard then $\nu(\t)$ is row semistandard. Finally, if $\T$ is a row
semistandard $\mu$--tableau of type $\nu$ then $\T^*$ is a row standard
$\mu$--tableau such that $\T=\nu(\T^*)$. Hence, $\E_{\mu\nu}$ is a complete set
of $(G_\mu,G_\nu)$--double coset representatives in~$\G$.

To complete the proof we need to show that if $\T\in \srssTab(\mu,\nu)$ then
$d_{\T^*}$ is an element of minimal length in the double coset
$G_\mu d_{\T^*}G_\nu$. For convenience, let $d=d_{|\T^*|}$. Then, 
$d\in\D_{\mup\nup}$ by (\ref{double Sym}). Now, by the last paragraph
$G_\mu d_{\T^*}G_\nu=\bigcup_\t G_\mu d_\t$, where $\t$ runs over the
row standard $\mu$--tableau $\t$ such that $\nu(\t)=\T$. By definition,
$d_{\T^*}=t_1^{\alpha_1}\dots t_n^{\alpha_n}d$, where if $x\in[\mup]$ then
$\T^*(x)=\zeta^{\alpha_i}i^d$ if and only if $\T(x)=\zeta^{\alpha_i}k$ and $i^d$ is in
row~$k$ of $\t^\nu$. Now suppose 
that $\t$ is any row standard $\mu$--tableaux such that $\nu(\t)=\T$. Then,
using (\ref{double Sym}) again,
$d_\t=t_1^{\beta_1}\dots t_n^{\beta_n}d_{|\t|}
     =t_1^{\beta_1}\dots t_n^{\beta_n}du$,
for some $u\in\Sym_\nu$ and where $\beta_i=\alpha_{i^w}$, for some
$w\in\Sym_n$ (since $\nu(\t)=\T$). Therefore, 
\begin{align*}
\ell(d_{\t})&=\beta_1+\dots+\beta_n+\ell(du)
             =\alpha_1+\dots+\alpha_n+\ell(du)\\
            &\ge \alpha_1+\dots+\alpha_n+\ell(d)=\ell(d_{\T^*}),
\end{align*}
with equality if and only if $u=1$. By Theorem~\ref{Dmu}, $d_\t$ is the
unique element of minimal length in the coset $G_\mu d_\t$, for each such
$\t$. Therefore, $\ell(d_{\T^*})\le\ell(g)$ for all $g\in G_\mu
d_{\T^*}G_\nu$ as claimed.
\end{proof}

Note that we are not claiming that each double coset of two reflection
subgroups of~$\G$ contains a unique element of minimal length. Indeed, the
proof of Proposition~\ref{double} shows that if~$\T$ is a row semistandard
$\mu$--tableau of type $\nu$ then the double coset $G_\mu d_{\T^*}G_\nu$
contains more than one element of minimal length if and only if there exist
integers $b,c$, not both zero, such that $m\zeta^b$ and~$m\zeta^c$ appear
in the same row of $\T$, for some $m\in\n$. For future comparison we make
this statement explicit.

If $d\in\D_{\mup\nup}$ let $\T_d\in\rssTab(\mup,\nup)$ be the unique row
semistandard tableau such that $d=d_{\T_d^*}$ as in (\ref{double Sym}). If
$X\subseteq\G$ let $X^{-1}=\set{g|g^{-1}\in X}$.

\begin{Lemma}\label{Emunu factorization I}
  Suppose that $\mu,\nu\in\SComp$. Then
  $$\E_{\mu\nu}=\coprod_{d\in\D_{\mup\nup}}\!\!
  \SET[55]t_1^{\alpha_1}\dots t_n^{\alpha_n}\in T_{-\mu\cap d(-\nu)}|%
     $\alpha_i\le \alpha_j$ whenever $i^d$ and  $j^d$ are in the same %
     row of $\T_d^*$ and the same row of $\t^\nu$|\,d
$$
Moreover, 
$$\E_\mu\cap\E_\nu^{-1}=\coprod_{d\in\D_{\mup\nup}}T_{-\mu\cap d(-\nu)}d
   =\set{e\in\G|\ell(e)\le\ell(g)\text{ for all }g\in G_\mu eG_\nu}$$
is the set of elements in~$\G$ which are of minimal length in their
$(G_\mu,G_\nu)$--double coset.
\end{Lemma}

\begin{proof}
Observe that $\D_{\mup\nup}=\D_{(-\mu)^+(-\nu)^+}$. Therefore, if
$d\in\D_{\mup\nup}$ then the signed composition $-\mu\cap d(-\nu)$ in the
statement of the Lemma makes sense by Definition~\ref{intersect scomp}.
(Note, however, that the two signed compositions $-\mu\cap d(-\nu)$ and
$-(\mu\cap d\nu)$ are \textit{not} equal in general.)

By Proposition~\ref{double}, we have
$\E_{\mu\nu}=\set{d_{\T^*}|\T\in\srssTab(\mu,\nu)}$. Fix a row
semistandard $\mu$--tableau~$\T$ of type~$\nu$. Then, as in the proof of
Proposition~\ref{double}, $d_{\T^*}=t_1^{\alpha_1}\dots t_n^{\alpha_n}d$, where
$d=d_{|\T^*|}\in\D_{\mup\nup}$ and, for all $x\in[\mup]$ if
$\T^*(x)=\zeta^{\alpha_i}i^d$ then $\T(x)=\zeta^{\alpha_i}k$ where $i^d$ is in
row~$k$ of $\t^\nu$. In particular, $\alpha_i=0$ if $t_i\in T_\mu$ or if
$t_{i^d}\in T_\nu$. Therefore, $\alpha_i>0$ only if 
$t_i\in T_{-\mu}\cap dT_{-\nu}d^{-1}=T_{-\mu\cap d(-\nu)}$. If
$t_i\in T_{-\mu\cap d(-\nu)}$ then the integer~$\alpha_i$ can take any value in
$\{0,\dots,r-1\}$ provided that this is compatible with $\T$ being row
semistandard. That is, we require that $\alpha_i\le \alpha_j$ whenever
$i^d$ and $j^d$ are in the same row of~$\T^*$ and in the same row of~$\t^\nu$.
This gives the decomposition of $\E_{\mu\nu}$ in the statement of Lemma.

For the final claim, suppose that $d\in\D_{\mup\nup}$ and let
$\T=\nu(\t^\mu d)$. By the last paragraph, if $t\in T$ then $\nu(\t^\mu
td_{\T^*})=\T$ if and only if $t\in T_{-\mu\cap d(-\nu)}$. By the last
paragraph again, if $t\in T_{-\mu\cap d(-\nu)}d$ then $td$ is an element of
minimal length in the double coset $G_\mu tdG_\nu$. That
$\E_\mu\cap\E_\nu^{-1}=\coprod_dT_{-\mu\cap d(-\nu)}d$ is now follows from
the definition of row semistandard tableaux.
\end{proof}

\begin{Cor}\label{Emunu factorization}
  Suppose that $\mu,\nu\in\SComp$ and $d\in\D_{\mup\nup}$. Then $\G$
  contains $|T_{-\mu\cap d(-\nu)}|$ elements of the form $t^\alpha d$ which
  are of minimal length in their $(G_\mu,G_\nu)$--double coset, for some
  $\alpha\in\Z_r^n$. Moreover, if $\T=\nu(\t^\mu d)$ then
  $|T_{-\mu\cap d(-\nu)}|=r^{\wt(\T)}$, where $\wt(\T)$ is the number of pairs
  $(i,j)$ such that $j$ appears in row~$i$ of~$\T$ and $\mu_i<0$ and
  $\nu_j<0$.
\end{Cor}

\begin{proof} That $|T_{-\mu\cap d(-\nu)}|$ counts the number of elements of the
  form $t^\alpha d$ which are of minimal length in their
  $(G_\mu,G_\nu)$--double coset is immediate from Lemma~\ref{Emunu
  factorization I}. The second claim 
  follows from the observation that the tableaux 
  $\set{\t^\mu td|t\in T_{-\mu\cap d(-\nu)}}$ differ only in that any of
  the numbers appearing in row~$i$ of $\t^\mu d$ and row~$j$ of~$\t^\nu$
  can be multiplied by arbitrary powers of~$\zeta$ whenever $\mu_i<0$ and
  $\nu_j<0$.
\end{proof}

\begin{Example}\label{row semistandard}
Suppose that $r\ge2$ and $n=5$ and let $\mu=(3,-2)$ and $\nu=(-2,-2,1)$.
Then the set of row semistandard $\mu$--tableaux $\T$ of type
$\nu$, together with the corresponding row standard tableau $\T^*$
and the coset representatives $d_{\T^*}\in\E_{\mu\nu}=\E_{(3,-2)(-2,-2,1)}$, 
is as follows (we set $d=d_{|\T^*|}$).
$$
\begin{array}{llccl}
\multicolumn1c{\T}&\multicolumn1c{\T^*}&\multicolumn1c{d_{\T^*}}
  &\multicolumn1c{|T_{-\mu\cap d(-\nu)}|}&\multicolumn1c{\mu\cap d\nu}\\\toprule
\Tab(1&1&2\cr2\zeta^a&3)&\Tab(1&2&3\cr4\zeta^a&5)
  &t_4^a&r&(-2,-1^3)\\[10pt]
\Tab(1&1&3\cr2\zeta^b&2\zeta^c)&\Tab(1&2&5\cr3\zeta^b&4\zeta^c)
  &t_4^bt_5^cs_3s_4&r^2&(-2,1,-2)\\[10pt]
\Tab(1&2&2\cr1\zeta^a&3)&\Tab(1&3&4\cr2\zeta^a&5)
  &t_4^as_3s_2&r&(-1,-2,-1^2)\\[10pt]
\Tab(1&2&3\cr1\zeta^a&2\zeta^b)&\Tab(1&3&5\cr2\zeta^a&4\zeta^b)
  &t_4^at_5^bs_3s_2s_4&r^2&(-1^2,1,-1^2)\\[10pt]
\Tab(2&2&3\cr1\zeta^b&1\zeta^c)&\Tab(3&4&5\cr1\zeta^b&2\zeta^c)
  &t_4^bt_5^cs_3s_2s_4s_3s_1s_2&r^2&(-2,1,-2)
\end{array}
$$
where $0\le a,b,c<r$ and $b\le c$.  We use exponentials in the signed
compositions to indicate consecutive repeated parts.  Therefore, there are
$2r^2+3r$ $(G_\mu,G_\nu)$--double cosets in~$\G$. When checking
the entries in this table observe that the signed composition $\mu\cap
d\nu=\mu\cap\nu d^{-1}$ can be computed by intersecting $G_\mu$ with the
`row stabilizer' of $\t^\nu d^{-1}$ as in the proof of 
Lemma~\ref{mud cap nu}. Note that $|T_{-\mu\cap d(-\nu)}|$ can be
computed without finding $-\mu\cap d(-\nu)$ by using 
Corollary~\ref{Emunu factorization}.
\end{Example}

\begin{Remark}
If $\mu$ and~$\nu$ are compositions of $n$ then
$\D_{\mu\nu}=\D_\mu\cap\D_\nu^{-1}$ is a complete set of minimal length
$(\Sym_\mu,\Sym_\nu)$--double coset representatives in~$\Sym_n$ by
(\ref{double Sym}). In contrast, it is not hard to show that
$\E_{\mu\nu}\subseteq\E_\mu\cap\E_\nu^{-1}$ with $\E_\mu\cap\E_\nu^{-1}$
being strictly bigger than $\E_{\mu\nu}$ in general. For example, if we
take $\mu=(3,-2)$ and $\nu=(-2^2,1)$ then $|\E_\mu\cap\E_\nu^{-1}|=3r^2+2r$,
whereas $|\E_{\mu\nu}|=2r^2+3r$ by Example~\ref{row semistandard}. So 
$\E_{\mu\nu}\subsetneq\E_\mu\cap\E_\nu^{-1}$ since $r>1$.
\end{Remark}

\section{The cyclotomic Solomon algebra}

Suppose that $R$ is a commutative ring (with one) and let $R\G$ be the
group ring of~$\G$ over~$R$. In this section we use the distinguished
coset representatives of the reflection subgroups of $\G$ to define an
analogue of Solomon's descent algebra for the complex reflection group
$\G$.

Recall that for each reflection subgroup $G_\mu$ of~$\G$ we have a
distinguished set $\E_\mu$ of right coset representatives, for
$\mu\in\SComp$. Define
$$E_\mu = \sum_{e\in\E_\mu} e\in R\G.$$
The main aim of this paper is to understand the subalgebra of $R\G$
which is generated by these elements.

\begin{Defn}\label{Solomon defn}
    Suppose that $r>1$. The \textbf{cyclotomic Solomon algebra} 
       $$\Sol(\G)=\Sol_R(\G)$$
    is the subalgebra of $R\G$ generated by $\set{E_\mu|\mu\in\SComp}$.
\end{Defn}

  From our definition, it is not clear what the dimension of $\Sol(\G)$ is
when $R$ is a field. In fact, we show in Theorem~\ref{cyc Solomon} below
that if $R$ is any ring then $\Sol(\G)$ is free as an $R$--module with
basis $\set{E_\mu|\mu\in\SComp}$.  We begin by taking advantage of the
factorization of $\E_\mu$ given by Theorem~\ref{Dmu}. To do this, for
$i=1,\dots,n$ and $\lambda\in\Comp$ define
$$ F_i=\sum_{k=0}^{r-1} t_i^k\qquad\text{and}\qquad
   D_\lambda=\sum_{d\in\D_\lambda}d,$$
Then $F_i$ and $D_\lambda$ are both elements of $R\G$.

\begin{Lemma}\label{F_is} 
Suppose that $1\le i,j\le n$ and that $w\in\Sym_n$. Then
\begin{enumerate}
\item $F_iF_j=F_jF_i$ and $F_i^2=rF_i$.
\item $F_iw=wF_{i^w}$.
\end{enumerate}
\end{Lemma}

\begin{proof} As $T$ is an abelian group part~(a) is true and
part~(b) is immediate from the definitions and~(\ref{T left}).
\end{proof}

Hence, if $1\le i\le n$ then $F_i$ is a multiple of an idempotent if the
characteristic of~$R$ does not divide~$r$ and, otherwise, it is a nilpotent
element of~$R\G$.

Suppose that $\mu\in\SComp$. In order to factorize $E_\mu$ set
\[F_{-\mu}=\prod_{t_i\in T_{-\mu}}F_i=\prod_{i:\mu_i<0}F_{\bar\mu_{i-1}+1}\cdots F_{\bar\mu_i}.\] 
Then, by Lemma~\ref{F_is}(a), $\(F_{-\mu})^2=r^{|\mu|^-}F_{-\mu}$. 

By Lemma~\ref{F_is}, $\Sym_n$ acts on
$\{F_1,\dots,F_n\}$ by conjugation. If $w\in\Sym_n$ and $i\in\n$ then 
we set $F_i^w=w^{-1}F_iw=F_{i^w}$. Similarly, if $\mu\in\SComp$ let
$$F_{-\mu}^w=\prod_{t_i\in T_{-\mu}}F_i^w.$$ 
Then $F_{-\mu}w=wF_{-\mu}^w$, for all $w\in\Sym_n$, by Lemma~\ref{F_is}(b).

\begin{Lemma}\label{D factorization}
    Suppose that $\mu\in\SComp$ is a signed composition of~$n$. Then:
  \begin{enumerate}
    \item $E_\mu=F_{-\mu} D_\mup$.
    \item If $w\in\Sym_{\mup}$ then $F_{-\mu}^ w=F_{-\mu}$, so that 
          $F_{-\mu} w=wF_{-\mu}$.
  \end{enumerate}
\end{Lemma}

\begin{proof} Part (a) is an immediate consequence of the factorization
$\E_\mu=T_{-\mu}\times\D_\mup$ of~$\E_\mu$ given by Theorem~\ref{Dmu}.
For part~(b), use Lemma~\ref{F_is}(b) and the fact that the elements of 
the two subgroups~$\Sym_\mu$ and~$T_{-\mu}$ commute.
\end{proof}

Definition~\ref{Solomon defn} is motivated by Solomon's~\cite{Solomon:descent}
definition of the descent algebra of a finite Coxeter group. As an
important special case, the \textbf{Solomon descent algebra}
$\Sol(\Sym_n)$ of~$\Sym_n$ is the subalgebra of $R\Sym_n$ generated by
$\set{D_\lambda|\lambda\in\Comp}$. The next result, due to Solomon,
shows that $\set{D_\lambda|\lambda\in\Comp}$ is  basis of $\Sol(\Sym_n)$.

\begin{Point}[\protect{Solomon~\cite[Theorem~1]{Solomon:descent}}]
\label{Solomon}\qquad\leavevmode\newline
\begin{enumerate}
\item\vskip-14pt The set $\set{D_\mu|\mu\in\Comp}$ is linearly
	independent in $\Sol(\Sym_n)$.
\item Suppose that $\mu$ and $\nu$ are composition of~$n$. Then
  $$D_\mu D_\nu = \sum_{d\in\D_{\mu\nu}}D_{\mu\cap d\nu}.$$
\end{enumerate}
\end{Point}

By the remarks before Lemma~\ref{mud cap nu}, part~(b) is equivalent to
the following formula:
$$D_\mu D_\nu=\sum_{\sigma\in\Comp}d_{\mu\nu\sigma}D_\sigma,$$
where $d_{\mu\nu\sigma}
    =\#\set{d\in\D_{\mu\nu}|\Sym_\sigma=\Sym_\mu \cap d^{-1}\Sym_\nu d}$. 
In fact, Solomon proved an analogous result for an arbitrary finite Coxeter
group~$W$, where the Young subgroups $\Sym_\mu$ are replaced with the
parabolic subgroups of~$W$ and $D_\mu$ by the sum of the `distinguished'
(right) coset representatives which are of minimal length in their coset.

As we now recall, part~(a) of Solomon's theorem is easy to prove. Recall
that $S=\{s_1,\dots,s_{n-1}\}$ and that if $w\in\Sym_n$ then $\Des(w)$ is
the descent set of $w$; see \ref{descent}. For each
composition $\mu\in\Comp$ let $S_\nu=\Pi_{-\nu}$, so that 
$S_\nu\subseteq S$. Now define $Y_\mu\in R\Sym_n\subset R\G$ by 
$$Y_\mu = \sum_{\substack{w\in\Sym_n\\\Des(w)=S_\mu}} w.$$ 
By definition, the descent sets partition $\Sym_n$, so the set 
$\set{Y_\mu|\mu\in\Comp}$ is linearly independent in $R\Sym_n$. 
By (\ref{Solomon}) again, we can write
$$D_\mu=\sum_{\substack{\nu\in\Comp\\S_\nu\subseteq S-S_\mu}} Y_\nu.$$
Hence, $\set{D_\mu|\mu\in\Comp}$ is a linearly independent subset of~$R\Sym_n$,
as claimed.

We build upon this idea to prove that the $E_\mu$'s are linearly
independent.

\begin{Prop}\label{basis theorem}
The set $\set{E_\mu|\mu\in\SComp}$ is linearly independent
in $\Sol(\G)$.
\end{Prop}

\begin{proof}
Suppose that there exist scalars $a_\mu\in R$ such that  
$$\sum_{\mu\in\SComp} a_\mu E_\mu =0.$$
By Lemma~\ref{D factorization}, $E_\mu = F_{-\mu} D_\mup$. Therefore, 
the last displayed equation becomes
\begin{align*}
0&= \sum_{\mu\in \SComp} a_\mu F_{-\mu} D_\mup
  = \sum_{\mu\in \SComp} a_\mu F_{-\mu}
          \sum_{\substack{\nu\in\Comp\\S_\nu\subseteq S-S_\mup}}Y_\nu\\
& = \sum_{\nu\in\Comp} \Big(
  \sum_{\substack{\mu\in\SComp\\S_\mup\subseteq S-S_\nu}} 
         a_\mu F_{-\mu}\,\Big)\, Y_\nu
\end{align*}
Now, $R\G=\bigoplus_{t\in T}tR\Sym_n$, as an $R$--module, and 
$\set{Y_\nu|\nu\in\Comp}$ is a linearly independent subset of $R\Sym_n$.
Therefore, for any composition $\nu\in\Comp$ we must have
$$
0=\sum_{\substack{\mu\in\SComp\\S_\mup\subseteq S-S_\nu}} a_\mu F_{-\mu}.
\leqno(*)$$
We use this equation to argue by induction on~$\nu$ to show that $a_\mu=0$
for all $\mu\in\SComp$. 

First suppose that $\nu=(n)$. Then $S_\nu=S$ and the summation in
$(*)$ becomes a sum over those signed compositions $\mu$ with
$S_\mup=\emptyset$. Hence, $\mup=(1^n)$ and $(*)$ becomes
$$0=\sum_{\substack{\mu\in\SComp\\\mup=(1^n)}}a_\mu F_{-\mu}
   =\sum_{\substack{\mu\in\SComp\\\Pi_\mu\subseteq\{t_1,\dots,t_n\}}}
        a_\mu F_{-\mu}.$$
Each monomial $t_{i_1}\dots t_{i_k}$, where $1\le i_1<\dots<i_k\le n$, 
occurs in a unique $F_{-\mu}$ when $\mup=(1^n)$.  Hence, $a_\mu=0$ for all
$\mu\in\SComp$ with $\mup=(1^n)$, as claimed.

Now suppose that $\nu\ne(n)$. By induction we may assume that $a_\mu=0$
whenever $S_\mup\subsetneq S-S_\nu$. Therefore, by $(*)$ we
have 
$$0=\sum_{\substack{\mu\in\SComp\\S_\mup=S-S_\nu}} a_\mu F_{-\mu}
   =\sum_{\substack{\mu\in\SComp\\\Pi_\mu-(S-S_\nu)\subseteq\{t_1,\dots,t_n\}}}
        a_\mu F_{-\mu}.$$
So, by exactly the same argument as before, $a_\mu=0$ whenever
$\mup=\nu$.  Hence, $a_\mu=0$, for all $\mu\in\SComp$, and
$\set{E_\mu|\mu\in\SComp}$ is linearly independent as required.  
\end{proof}

The next result that we need amounts to a proof of part~(b) of Solomon's
theorem~(\ref{Solomon}). Once again, we state the result only for the
symmetric group even though it is valid for an arbitrary finite Coxeter
group. All of the results quoted in (\ref{BBHT}) follow easily from the
fact that $\Sym_{\mu d\cap \nu}=d^{-1}\Sym_\mu d\cap\Sym_\nu$, for
$d\in\D_{\mu\nu}$.

If $\mu,\nu\in\Comp$ and $\Sym_\nu\subseteq\Sym_\mu$ then we write
$\nu\subseteq\mu$ and set $\D_\nu^\mu=\D_\nu\cap\Sym_\mu$. It is easy to
check that $\D_\nu^\mu$ is a complete set of coset representatives for
$\Sym_\nu$ in $\Sym_\mu$.

\begin{Point}[\protect{Bergeron, Bergeron, Howlett and
    Taylor~\cite[Lemmas 2.2 and 2.4]{BBHT}}]\label{BBHT}
Suppose that $\mu$ and $\nu$ are compositions of $n$. Then
\begin{enumerate}
  \item If $\sigma\subseteq\nu$ then $\D_\sigma=\D_\sigma^\nu\D_\nu$.
  \item $\D_\mu=\coprod_{d\in\D_{\mu\nu}}d\D_{\mu d \cap\nu}^\nu$.
  \item If $d\in\D_\mu$ and $\mu d$ is a composition
  of~$n$ $($that is, $d^{-1}\Sym_\mu d=\Sym_\sigma$ for some
  $\sigma\in\Comp)$, then $\D_\mu=d\D_{\mu d}$.
\end{enumerate}
\end{Point}

We can now establish one of the main results of this paper. 

\begin{Thm}\label{cyc Solomon}
Suppose that $r>1$ and that $\mu$ and $\nu$ are signed compositions of~$n$. Then 
$$E_\mu E_\nu
    =\sum_{d\in\D_{\mup\nup}}|T_{-\mu\cap d(-\nu )}|\,E_{\mu\cap d\nu}.
$$
\end{Thm}

\begin{proof} We use most of the results in this section to compute 
  $E_\mu E_\nu$:
\begin{align*} E_\mu E_\nu 
  &= F_{-\mu} D_\mup F_{-\nu} D_\nup,
                &\text{by Lemma~\ref{D factorization}(a)},\\
  &= \sum_{d\in\D_{\mup\nup}} 
        F_{-\mu} dD^\nup_{\mup d\cap\nup}F_{-\nu} D_\nup,
                &\text{by Lemma~\ref{BBHT}(b),}\\
  &= \sum_{d\in\D_{\mup\nup}} 
        F_{-\mu} dF_{-\nu} D^\nup_{\mup d\cap\nup}D_\nup,
                &\text{by Lemma~\ref{D factorization}(b)},\\
  &= \sum_{d\in\D_{\mup\nup}} F_{-\mu} F_{-\nu}^{d^{-1}} dD_{\mup d\cap\nup},
                &\text{by Lemma~\ref{BBHT}(b)},\\
  &= \sum_{d\in\D_{\mup\nup}} F_{-\mu} F_{-\nu}^{d^{-1}} D_{\mup\cap d\nup},
                &\text{by Lemma~\ref{BBHT}(c)}.\\
\end{align*}
Fix $d\in\D_{\mup\nup}$ and consider $F_{-\mu}F_{-\nu}^{d^{-1}}$. Now
$F_i^2=rF_i=|T_i|F_i$, for $1\le i\le n$. So,
$$F_{-\mu}F_{-\nu}^{d^{-1}}=|T_{-\mu}\cap dT_{-\nu}d^{-1}| 
\prod_{t_i\in T_{-\mu}\cap dT_{-\nu}d^{-1}}F_i.$$
First,
$T_{-\mu}\cap dT_{-\nu}d^{-1}=T_{-\mu\cap d(-\nu)}$
since $d\in\D_{\mup\nup}=\D_{(-\mu)^+(-\nu)^+}$.
Next, the subgroup of~$T$ generated by $T_{-\mu}$ and 
$dT_{-\nu}d^{-1}$ is $T_{-(\mu\cap d\nu)}$ since 
$t_i\in T_{-(\mu\cap d\nu)}\cong T/T_{\mu\cap d\nu}$
if and only if $t_i\notin T_\mu$ and $t_i\notin dT_\nu d^{-1}$.
Therefore, 
$F_{-\mu}F_{-\nu}^{d^{-1}}=|T_{-\mu\cap d(-\nu)}|F_{-(\mu\cap d\nu)}$.
Hence, using Lemma~\ref{D factorization} once more,
$$E_\mu E_\nu
    =\sum_{d\in\D_{\mup\nup}}|T_{-\mu\cap d(-\nu)}|\,
            F_{-(\mu\cap d\nu)}D_{(\mu\cap d\nu)^+}
    =\sum_{d\in\D_{\mup\nup}}|T_{-\mu\cap d(-\nu)}|\,E_{\mu\cap d\nu},
$$
as required.
\end{proof}

Corollary~\ref{Emunu factorization} shows that $|T_{-\mu\cap d(-\nu)}|$ is
equal to the number of elements of minimal length in the double cosets of
the form $G_\mu t^\alpha dG_\nu$, for $\alpha\in\Z_n^r$. This gives a
combinatorial interpretation of the structure constants of~$\Sol(\G)$ and
shows that Theorem~\ref{cyc Solomon} a direct generalization of
(\ref{Solomon})(b). A second combinatorial interpretation of the integers
$|T_{-\mu\cap d(-\nu)}|$ is given in Proposition~\ref{combII}
below.

Combining Theorem~\ref{cyc Solomon} and Proposition~\ref{basis theorem} we
obtain the following.

\begin{Cor}\label{dimension}
Suppose that $r>1$. The cyclotomic Solomon algebra  $\Sol(\G)$ is a
subalgebra of $R\G$ which is free as an $R$--module of rank 
$2\cdot 3^{n-1}$.
\end{Cor}

\begin{Example}\label{multiplication example 1}
Suppose that $r>1$. Then, by Example~\ref{row semistandard}, we have
$$E_{(3,-2)}E_{(-2^2,1)}
  =2r^2E_{(-2,1,-2)}+rE_{(-2,-1^3)}+rE_{(-1,-2,-1^2)}+r^2E_{(-1^2,1,-1^2)}.$$
See Example~\ref{multiplication example 2} for a second way of computing
this product using Proposition~\ref{combII}.

Notice that by (\ref{Solomon}) and Theorem~\ref{cyc Solomon} we can recover
the multiplication in $\Sol(\Sym_n)$ by setting $r=1$ and identifying $\mu$
and $\mup$, for all $\mu\in\SComp$, so that
$$D_{(3,2)}D_{(2^2,1)}=2D_{(2,1,2)}+D_{(2,1^3)}+D_{(1,2,1^2)}+D_{(1^5)}.$$
\end{Example}

\section{The generic cyclotomic Solomon algebra}

By Theorem~\ref{cyc Solomon}, if $r>1$ then the structure constants of the
algebra $\Sol(\G)$ are polynomials in $r$. Consequently, the algebras
$\set{\Sol(\G)|r\ge2}$ admit a simultaneous deformation (while $n$ is fixed).

Recall from Corollary~\ref{Emunu factorization} that if $r\ge2$ and
$\mu,\nu\in\SComp$ then 
$|T_{-\mu\cap d(-\nu)}|=r^{\wt(\T_d)}$, where $\T_d=\nu(\t^\mu d)$. 

Let $x$ be an indeterminate over $\Z$ and suppose that
$\mu,\nu,\sigma\in\SComp$.  Define polynomials 
$d_{\mu\nu\sigma}(x)\in\N[x]$ by
$$d_{\mu\nu\sigma}(x)=\sum_{\substack{d\in\D_{\mup\nup}\\\sigma=\mu\cap d\nu}}
x^{\wt(\T_d)}.$$

We abuse notation and consider $d_{\mu\nu\sigma}(x)$ to be a
polynomial over~$R$. For $q\in R$ we let $d_{\mu\nu\sigma}(q)$ be the
evaluation of this polynomial at $q$.
Then, by Theorem~\ref{cyc Solomon},
$$E_\mu E_\nu=\sum_{\sigma\in\SComp} d_{\mu\nu\sigma}(r) E_\sigma.$$

\begin{Defn}\label{generic}
  Suppose that $n\ge1$ and that $R$ is a commutative ring. The 
  \textbf{cyclotomic Solomon algebra with parameter $q\in R$} is the
  $R$--algebra  $\Sol_q(n)=\Sol_{R,q}(n)$ with generating set
  $\set{E_\mu|\mu\in\SComp}$ and relations
  $$E_\mu E_\nu=\sum_{\sigma\in\SComp}d_{\mu\nu\sigma}(q)E_\sigma,$$
  for $\mu,\nu\in\SComp$. The \textbf{generic cyclotomic Solomon algebra}
  is the $\Z[x]$--algebra $\Sol_x(n)$, where $x$ is an indeterminate 
  over~$\Z$.
\end{Defn}

We are abusing notation slightly in Definition~\ref{generic} because from
here onwards $E_\mu$ is a generator of~$\Sol_q(n)$ and not necessarily the
element defined in the previous section. This abuse is justified by the
following result.

\begin{Cor}\label{specialize}
  Suppose that $q=r\cdot 1_R$, where $r>1$. Then $\Sol_q(n)$ and $\Sol(\G)$ are
  canonically isomorphic $R$--algebras where the isomorphism 
  $\Sol_q(n)\to\Sol(\G)$ is given by
  $E_\mu\mapsto E_\mu$, for $\mu\in\SComp$.
\end{Cor}

\begin{proof}
By Theorem~\ref{cyc Solomon} there is a natural surjection
$\Sol_q(n)\longrightarrow\Sol(\G)$. By Corollary~\ref{dimension} 
this map is an isomorphism.
\end{proof}

The explicit description of the algebra $\Sol(\G)$ as a subalgebra of the
group algebra $R\G$ makes the algebra $\Sol(\G)$ slightly easier to work
with than the more general algebras $\Sol_q(n)$. For example, we know that
$E_\mu=F_{-\mu} D_\mup$ in $R\G$ but we have no such factorization in
general. As we will soon see, however, almost all of the properties of the
algebras $\Sol(\G)$ hold for the algebras $\Sol_q(n)$.

\begin{Prop}\label{parameters}
  Suppose that $n\ge1$ and that $q\in R$. Then
  \begin{enumerate}
    \item $\Sol_q(n)$ is free as an $R$--module with basis
      $\set{E_\mu|\mu\in\SComp}$.
      In particular, $\Sol_q(n)$ has rank $2\cdot3^{n-1}$.
    \item $\Sol_q(n)\cong\Sol_x(n)\otimes_{\Z[x]}R$, where $R$ is
      considered as a $\Z[x]$--module by letting~$x$ act on~$R$ as 
      multiplication by $q$ $($and $1\in\Z$ acts as multiplication by $1_R)$.
    \item $\Sol_q(n)$ is a unital associative $R$--algebra with multiplicative
      identity $E_{(n)}$. 
  \end{enumerate}
\end{Prop}

\begin{proof}
First consider the generic Solomon algebra over $\Z[x]$. Suppose that 
$$\sum_{\mu\in\SComp}f_\mu(x)E_\mu=0,$$
for some $f_\mu(x)\in R[x]$. Then $f_\mu(r)=0$, for $r=2,3,4\dots$ and all
$\mu\in\SComp$, by Corollary~\ref{specialize} and 
Proposition~\ref{basis theorem}. As non--zero polynomials have only
finitely many roots, we conclude that $f_\mu(x)=0$, for all $\mu\in\SComp$.
Consequently, $\Sol_x(n)$ is free as a $\Z[x]$--module with basis
$\set{E_\mu|\mu\in\SComp}$. 

Now fix $q\in R$ and consider $R$ as a $\Z[x]$--module by letting $x$ act
on $R$ as multiplication by~$q$ (and $1\in\Z$ act as multiplication
by~$1_R$). Then the $R$--algebra $\Sol_x(n)\otimes_{\Z[x]}R$ is free as an
$R$--module with basis $\set{E_\mu\otimes1|\mu\in\SComp}$ and it satisfies
the relations of $\Sol_q(n)$. As $\Sol_q(n)$ is spanned by the elements
$\set{E_\mu|\mu\in\SComp}\subseteq\Sol_q(n)$ it follows that
$\Sol_q(n)\cong\Sol_x(n)\otimes_{\Z[x]}R$. This proves (a) and (b).
  
To prove (c) it is now enough to prove the corresponding statements for the
generic Solomon algebra $\Sol_x(n)$. We first show that $E_{(n)}$ is the
identity element of $\Sol_x(n)$. This is equivalent to the polynomial
identities 
$$d_{\mu(-n)\alpha}(x)=\delta_{\mu\alpha}=d_{\alpha(-n)\mu}(x),$$
for all $\mu,\alpha\in\SComp$. All of these identities follow directly 
from the definitions because $T_{(-n)\cap d(-\nu)}=1=T_{-\mu\cap d'(-n)}$, 
for all $\mu,\nu\in\SComp$, $d\in\D_{(n)\mup}$ and $d'\in\D_{\nup(n)}$. 
Similarly, the associativity of 
$\Sol_x(n)$ is equivalent to the polynomial identities
$$\sum_{\alpha,\beta\in\SComp}d_{\mu\nu\alpha}(x)d_{\alpha\sigma\beta}(x)
   =\sum_{\alpha,\beta\in\SComp}d_{\mu\alpha\beta}(x)d_{\nu\sigma\alpha}(x),$$
for all $\mu,\nu,\sigma\in\SComp$. As in the first paragraph of the
proof, by Corollary~\ref{specialize} these identities hold when
$x=2,3,\dots$since the algebras $\Sol(\G)$ are associative for $r\ge2$. As
these identities hold for infinitely many values of~$x$, they lift to the
required polynomial identities.
\end{proof}

Part~(b) of the Proposition justifies our calling the $\Z[x]$--algebra
$\Sol_x(n)$ the \textit{generic} cyclotomic Solomon algebra.

As we next describe, the algebras $\Sol_q(n)$ have many interesting subalgebras. 

\begin{Lemma}\label{refinement}
  Suppose that $d_{\mu\nu\sigma}(q)\ne0$, for $\mu,\nu,\sigma\in\SComp$. Then 
  $\Pi_\sigma=\Pi_\mu\cap d\Pi_\nu d^{-1}$, for some $d\in\D_{\mup\nup}$
\end{Lemma}

\begin{proof} 
By definition, the polynomial $d_{\mu\nu\sigma}(x)$ is non--zero only if
$G_\sigma=G_\mu\cap dG_\nu d^{-1}$ for some $d\in\D_{\mup\nup}$. Consequently,
if $d_{\mu\nu\sigma}(q)\ne0$ then
$\Pi_\sigma=\Pi_\mu\cap d\Pi_\nu d^{-1}$, for some
$d\in\D_{\mup\nup}$.
\end{proof}

Notice, in particular, that this implies that the poset structure on
$\SComp$ given by defining $\mu\preceq\nu$ whenever $\Pi_\nu\subseteq
\Pi_\mu$ is compatible with the ideal structure of $\Sol_q(n)$.

\begin{Prop}\label{filtration}
Suppose that $n\ge1$ and that $q\in R$. Then $\Sol_q(n)$ has a 
filtration by two--sided ideals
  $$\Sol_q(n)=\S_0\supset\dots\supset \S_n\supset0$$
where $\S_i$ is the $R$--submodule of $\Sol_q(n)$ with basis
$\set{E_\mu|\mu\in\SComp\text{ such that }|\mu|^+\ge i}$, for $i=0,\dots,n$. 
\end{Prop}

\begin{proof}By Lemma~\ref{refinement}, $d_{\mu\nu\sigma}(q)\ne0$ only if
  $\Pi_\sigma=\Pi_\mu\cap d\Pi_\nu d^{-1}$, for some $d\in\D_{\mup\nup}$. 
  Consulting the definitions, $|\mu|^+=|\Pi_\mu\cap T|$. Therefore, 
  $d_{\mu\nu\sigma}(q)\ne0$ only if $|\sigma|^+\le\min\{|\mu|^+,|\nu|^+\}$. 
  Hence, $\S_i$ is a two--sided ideal of 
  $\Sol(\G)$, for $0\le i\le n$, and the Proposition follows.
\end{proof}

\begin{Prop}\label{subalgebras}
Suppose that $n\ge1$ and that $q\in R$. Let 
\begin{align*}
  \Sol_q^+(n)&=\sum_{\mu\in\Comp}RE_\mu\\
  \Sol_q^\pm(n)&=\sum_{\pm\mu\in\Comp}RE_\mu\\
  \Sol_q^1(n)&=\sum_{\substack{\mu\in\SComp\\\mu_i>0\text{ for }i>1}}RE_\mu.
\end{align*}
Then $\Sol_q^+(n)$, $\Sol_q^\pm(n)$ and $\Sol^1_q(n)$ are all subalgebras
of $\Sol_q(n)$. Moreover, $\Sol_q^+(n)$ is naturally isomorphic to
$\Sol(\Sym_n)$ via the $R$--linear map $E_\mu\mapsto D_\mu$, for
$\mu\in\Comp$.
\end{Prop}

\begin{proof}All of these results can be proved directly using the
  definition of the polynomials $d_{\mu\nu\sigma}(x)$, for
  $\mu,\nu,\sigma\in\SComp$. Note that $\Sol_q^+(n)=\S_n$ in the notation
  of Proposition~\ref{filtration}, so in this case the result is already
  known. The isomorphism $\Sol_q^+(n)\cong\Sol(\Sym_n)$ is trivial because
  if $\mu\in\Comp$ then $T_{-\mu}=1$, so that
  $E_\mu=D_\mu$ by Lemma~\ref{D factorization}.
\end{proof}

\section{The representation theory of $\Sol_q(n)$} 
In this section we construct all of the irreducible representations of the
algebras $\Sol_q(n)$ over an arbitrary field. Even though $\Sol_q(n)$ is,
in general, not commutative, it turns out that every irreducible
$\Sol_q(n)$--module is one dimensional --- so that $\Sol_q(n)$ is a basic
algebra for all $n$ and $q$. As an application of these results we give a
basis for the radical of $\Sol_q(n)$ when $R$ is an arbitrary field.

Let $\sim$ be the equivalence relation on the set of signed compositions
where two signed compositions are $\sim$--equivalent if one can be obtained
by reordering the parts of the other. More explicitly, if
$\lambda=(\lambda_1,\dots,\lambda_k)$ and $\mu=(\mu_1,\dots,\mu_l)$ then
$\lambda\sim\mu$ if and only if $k=l$ and $\lambda_i=\mu_{i^v}$, for some
$v\in\Sym_k$. 

\begin{Lemma}\label{conjugate}
  Suppose that $\lambda,\mu\in\SComp$. Then the following are equivalent:
  \begin{enumerate}
    \item $\lambda\sim\mu$;
    \item $G_\lambda=w^{-1}G_\mu w$, for some $w\in\Sym_n$;
    \item $G_\lambda=g^{-1}G_\mu g$, for some $g\in\G$.
  \end{enumerate}
\end{Lemma}

\begin{proof}
  We leave the proof for the reader.
\end{proof}

\begin{Lemma}\label{equivalence}
Suppose that $\lambda,\mu\in\SComp$. Then
\begin{enumerate}
\item If $\mu\not\sim\lambda$ then 
  $d_{\mu\lambda\mu}(q)\ne0$ only if $|\Pi_\lambda|>|\Pi_\mu|$.
\item If $\mu\sim\lambda$ then $d_{\mu\alpha\mu}(q)=d_{\lambda\alpha\lambda}(q)$,
  for all $\alpha\in\SComp$.
\end{enumerate}
\end{Lemma}

\begin{proof}
By Lemma~\ref{refinement} $d_{\mu\lambda\mu}(q)\ne0$ only if
$\Pi_\mu=\Pi_\mu\cap d\Pi_\lambda d^{-1}$, for some
$d\in\D_{\mup\lamp}$.  Hence, part~(a) follows since $\lambda\not\sim\mu$.

Consulting the definition of the polynomials $d_{\mu\nu\sigma}(x)$, 
to prove (b) it is enough to show that if $r\ge2$ then in the group $\G$
we have
$$\sum_{\substack{d\in\D_{\lamp\alpp}\\\lambda=\lambda\cap d\alpha}}
     |T_{-\lambda\cap d(-\alpha)}|
    =\sum_{\substack{d\in\D_{\mup\alpp}\\\mu=\mu\cap d\alpha}}
     |T_{-\mu\cap d(-\alpha)}|.\leqno(\dag)$$
We prove
this by showing that the `obvious' bijection
$\D_{\lamp\alpp}\longrightarrow\D_{\mup\alpp}$ preserves each of the
summands in this equation.

First note that by Lemma~\ref{conjugate} we can find an element
$w\in\Sym_n$ such that $G_\lambda=w^{-1}G_\mu w$ since $\lambda\sim\mu$. That is,
$T_\lambda\Sym_{\lamp}=w^{-1}T_\mu w\cdot w^{-1}\Sym_{\mup} w$, so that
$T_\lambda=w^{-1}T_\mu w$ and $\Sym_{\lamp}=w^{-1}\Sym_{\mup} w$.
Consequently,  the map
$\Sym_{\lamp}\backslash\Sym_n\slash\Sym_{\alpp}
       \to\Sym_{\mup}\backslash\Sym_n\slash\Sym_{\alpp}$ given by
$C\mapsto wC$ defines a bijection since if
$d\in\D_{\lamp\alpp}$ then 
$\Sym_{\lamp} d\Sym_{\alpp}=w^{-1}\Sym_{\mup} wd\Sym_{\alpp}$. Let $d\mapsto d'$ be map
from $\D_{\lamp\alpp}$ to $\D_{\mup\alpp}$ determined 
by
$\Sym_{\mup} wd\Sym_{\alpp}=\Sym_{\mup} d'\Sym_{\alpp}$.

Now fix $d\in\D_{\lamp\alpp}$ such that $\lambda=\lambda\cap d\alpha$. Then
\begin{align*}
  T_{-\lambda\cap d(-\alpha)}
      &=T_{-\lambda}\cap dT_{-\alpha}d^{-1}\\
      &=w^{-1}T_{-\mu}w\cap dT_{-\alpha}d^{-1}\\
      &=w^{-1}\Big(T_{-\mu}\cap wdT_{-\alpha}(wd)^{-1}\Big)w.\\
\intertext{Write $wd=ud'v$, for $u\in\Sym_{\mup}$ and $v\in\Sym_{\alpp}$. Then
we have}
  T_{-\lambda\cap d(-\alpha)}
  &=w^{-1}\Big(T_{-\mu}\cap (ud'v)T_{-\alpha} v^{-1}(d')^{-1}u^{-1}\Big)w,\\
  &=w^{-1}u\Big(T_{-\mu}\cap d'T_{-\alpha}(d')^{-1}\Big)u^{-1}w,\\
  &=w^{-1}u\Big(T_{-\mu}\cap T_{d'(-\alpha)}\Big)u^{-1}w,\\
\end{align*}
where the second equality follows because $\Sym_\mu$ normalizes
$T_{-\mu}$ and the last equality follows because $\Sym_\alpha$ 
normalizes~$T_{-\alpha}$. Hence, we have
shown that $|T_{-\lambda\cap d(-\alpha)}|=|T_{-\mu\cap d'(-\alpha)}|$,
for all $d\in\D_{\lamp\alpp}$. This establishes $(\dag)$, so the
Lemma is proved.
\end{proof}

\begin{Thm}\label{irreducibles}
  Suppose that $R$ is a field, $q\in R$ and $n\ge0$.
  \begin{enumerate}
    \item If $\lambda\in\SComp$ then $\Sol_q(n)$ has a unique
      one dimensional representation $I(\lambda)$ upon which $E_\alpha$
      acts as multiplication by $d_{\lambda\alpha\lambda}(q)$, for
      $\alpha\in\SComp$.
    \item Every irreducible representation of $\Sol_q(n)$ is isomorphic to
      $I(\lambda)$, for some $\lambda\in\SComp$.
    \item If $\lambda\sim\mu$ then $I(\lambda)\cong I(\mu)$.
  \end{enumerate}
\end{Thm}

\begin{proof}
Choose a total order $\ge$ on $\SComp$ such that $|\Pi_\lambda|\ge|\Pi_\mu|$
whenever $\lambda>\mu$, for $\lambda,\mu\in\SComp$.
Let $\S_\lambda$ be the $R$--submodule of $\Sol_q(n)$ with basis
$\set{E_\mu|\lambda\ge\mu\in\SComp}$ and let $\S_\lambda'$ be the
$R$--submodule with basis $\set{E_\mu|\lambda>\mu\in\SComp}$. Then
$\S_\lambda$ and $\S_\lambda'$ are both right $\Sol_q(n)$--modules by
Lemma~\ref{refinement}. Hence the quotient module
$I(\lambda)=\S_\lambda/\S_\lambda'=R(E_\lambda+\S_\lambda')$ is one
dimensional $\Sol_q(n)$--module. By definition, if $\alpha\in\SComp$ then
$E_\alpha$ acts on $I(\lambda)$ as multiplication by
$d_{\lambda\alpha\lambda}(q)$. Hence, $I(\lambda)$ is the one dimensional
$\Sol_q(n)$--module described in part~(a).  

Now suppose that $\SComp=\{\lambda_1>\lambda_1>\dots>\lambda_N\}$, where
$N=2\cdot3^{n-1}=\dim\Sol_q(n)$. Then
$$\Sol_q(n)=\S_{\lambda_1}\supset\S_{\lambda_2}\supset\dots\supset
      \S_{\lambda_N}\supset0$$
is a filtration of $\Sol_q(n)$ by two--sided ideals with quotients
$\S_{\lambda_i}/\S_{\lambda_{i+1}}\cong I(\lambda_i)$, since
$\S_{\lambda_{i+1}}=\S_{\lambda_i}'$. As every irreducible $\Sol_q(n)$--module
arises as a composition factor of $\Sol_q(n)$ part~(b) now follows.

Finally, if $\lambda\sim\mu$ then $I(\lambda)\cong I(\mu)$ by 
Lemma~\ref{equivalence}(b). Hence, part~(c) holds.
\end{proof}

\begin{Cor}Every field is a splitting field for $\Sol_q(n)$.
\end{Cor}

\begin{proof}Suppose that $D$ is an irreducible $\Sol_q(n)$--module. Then
  $D$ is one dimensional by the Proposition, and hence absolutely
  irreducible.
\end{proof}

If $A$ is an algebra over a field then let $\Rad A$ be its
\textbf{radical}. Thus, $\Rad A$ is the unique maximal nilpotent ideal of
$A$ and $A$ is semisimple if and only if $\Rad A=0$. Recall that $a\in A$
is \textbf{nilpotent} if $a^k=0$, for some $k>0$, whereas an ideal~$I$
of~$A$ is \textbf{nilpotent} if~$I^k=0$ for some $k>0$.

\begin{Cor}\label{nilpotent}
  Suppose that $R$ is a field. Then $\Rad\Sol_q(n)$ is the set of
  nilpotent elements in~$\Sol_q(n)$.
\end{Cor}

\begin{proof}
  By definition every element of $\Rad\Sol_q(n)$ is nilpotent. To prove the
  converse let $M$ be the number of irreducible $\Sol_q(n)$--modules. By
  Theorem~\ref{irreducibles} every irreducible $\Sol_q(n)$--module is one
  dimensional. Therefore, $\Sol_q(n)/\Rad\Sol_q(n)\cong R^M$ by the
  Wedderburn Theorem. In particular, $\Sol_q(n)/\Rad\Sol_q(n)$ contains no
  nilpotent elements, so the result follows.
\end{proof}

\begin{Cor}
  Suppose that $R$ is a field and $q\in R$. Then $\Sol_q(n)$ is semisimple
  if and only if $n=1$ and $q\ne0$.
\end{Cor}

\begin{proof}
If $n\ge2$ then $\Sol_q(n)$ is not semisimple because there exist
distinct signed compositions $\lambda,\nu\in\SComp$ such that
$\lambda\sim\mu$. Therefore, $E_\lambda-E_\mu\in\Rad\Sol_q(n)$, so that
$\Rad\Sol_q(n)\ne0$. If $n=1$ then a quick calculation verifies that
$I(1)\cong I(-1)$ if and only if $q=0$ which implies the result.
\end{proof}

Each $\sim$--equivalence class of $\SComp$ contains a unique signed
composition $\mu=(\mu_1,\dots,\mu_k)$ such that $\mu_1\ge\dots\ge\mu_k$.
If $\mu\in\SComp$ and $\mu_1\ge\dots\ge\mu_k$ then we call $\mu$ a
\textbf{signed partition} of~$n$.  Let $\SPart$ be the set of all signed
partitions of~$n$. By the remarks above, the $\Sym_n$--conjugacy classes of
reflection subgroups of $\G$ are indexed by the signed partitions of~$n$.
We note that $\SPart$ is naturally in bijection with the set of
bipartitions of~$n$, however, for us the signed partitions are more natural
because  we have already defined a reflection subgroup $G_\lambda$ for each
$\lambda\in\SPart$.

\begin{Thm}\label{char zero irreds}
Suppose that $R$ is a field of characteristic zero and that $q\ne0$. Then 
$$\set{I(\lambda)|\lambda\in\SPart}$$
is a complete set
of pairwise non--isomorphic irreducible $\Sol_q(n)$--modules.
\end{Thm}

\begin{proof}
As the $\sim$--equivalence classes of $\SComp$ are indexed by the
signed partitions of~$n$, $\set{I(\lambda)|\lambda\in\SPart}$ is a
complete set of irreducible $\Sol_q(n)$--modules by parts~(b) and~(c)
of Theorem~\ref{irreducibles}. It remains then to show that if
$\lambda,\mu\in\SPart$ then $I(\lambda)\not\cong I(\mu)$ if
$\lambda\ne\mu$. Now,~$R$ is a field  of characteristic zero and $q\ne0$, so
$d_{\lambda\nu\lambda}(q)\ne0$ if and only if
$d_{\lambda\nu\lambda}(x)\ne0$, for $\lambda,\nu\in\SComp$. However,
$d_{\lambda\lambda\lambda}(x)\in 1+x\N[x]$ since
$1\in\D_{\lamp\lamp}$ and
$\Pi_\lambda=\Pi_\lambda\cap 1\cdot\Pi_\lambda\cdot 1^{-1}$. Therefore,
$d_{\lambda\lambda\lambda}(q)\ne0$ and so, using
Lemma~\ref{equivalence}(a) again, if
$\lambda\ne\mu$ then $I(\lambda)\not\cong I(\mu)$.
\end{proof}

\begin{Cor}\label{radical}
Suppose that $R$ is a field of characteristic zero and $q\ne0$. Then 
$$\set{E_\lambda-E_\mu|\lambda\in\SPart, \mu\in\SComp, \lambda\sim\mu
         \text{ and }\lambda\ne\mu}$$
is a basis of $\Rad\Sol_q(n)$. Consequently,
  $\dim\Sol_q(n)/\Rad\Sol_q(n)=|\SPart|$.
\end{Cor}

\begin{proof}
Suppose that  $\lambda\sim\mu$ where $\lambda\in\SPart$, $\mu\in\SComp$ and
$\lambda\ne\mu$. Then, by Theorem~\ref{char zero irreds} and
Lemma~\ref{equivalence}, $E_\lambda-E_\mu$ acts as multiplication by zero
on every irreducible $\Sol_q(n)$--module. Therefore, $E_\lambda-E_\mu$
belongs to $\Rad\Sol_q(n)$ whenever $\lambda\sim\mu$. Consequently,
$\dim\Sol_q(n)/\Rad\Sol_q(n)\le|\SPart|$. However, 
$\dim\Sol_q(n)/\Rad\Sol_q(n)=|\SPart|$ by Theorem~\ref{char zero irreds}, 
so the result follows.
\end{proof}

Suppose that $R$ is a field of characteristic zero and that $q\ne0$. 
Define the \textbf{character table} of~$\Sol_q(n)$ to be the matrix 
$$\mathbf C_q(n)=\(d_{\lambda\mu\lambda}(q)\)_{\lambda,\mu\in\SPart}.$$
Then $\mathbf C_q(n)$ is the character table of $\Sol_q(n)/\Rad\Sol_q(n)$, by
Theorem~\ref{char zero irreds}, so it completely determines the
maximal semisimple quotient of $\Sol_q(n)$.  The character table $\mathbf C_q(n)$ is
explicitly known for all $q\ne0$ and all $n\ge1$ since the polynomials
$d_{\lambda\mu\sigma}(x)$ are explicitly known for all
$\lambda,\mu,\sigma\in\SComp$ by Corollary~\ref{Emunu factorization}.

\begin{Example}\label{character}
Suppose that $R$ is a field of characteristic zero and that $q=2=n$. Then 
$\Sol_2(2)\cong\Sol(G_{2,2})$ and
the character table $\mathbf C_2(2)$ of $\Sol_2(2)$ is the following matrix.
$$\begin{array}{l|*{6}{c}}
  &(2) & (1^2) & (1,-1) & (-2) & (-1^2)\\\toprule
(2)   &1&&&&&\\
(1^2) &1&2&&&&\\
(1,-1)&1&2&2&&&\\
(-2)  &1&.&.&4&&\\
(-1^2)&1&2&4&4&8\\
\end{array}$$
As all of the diagonal entries of $\mathbf C_q(2)$ are powers
of~$2$ it follows that if $R$ is any field of characteristic
\textit{different from}~$2$ then
$\set{I(\lambda)|\lambda\in\SPart}$ is a complete set of pairwise
non--isomorphic irreducible $\Sol_q(2)$--modules. If $R$ is a field of
characteristic~$2$ then $I(2)$ is the only irreducible $\Sol_q(2)$--module.
This is in agreement with Theorem~\ref{general irreducibles} below.

By comparing the character table of $\Sol(G_{2,2})$ with the character table
of the group $G_{2,2}$ (the Coxeter group of type $B_2$) it is easy to see
that there cannot be a ring homomorphism from $\Sol(G_{2,2})$ into the
character ring of $G_{2,2}$. This is in marked contrast with the Solomon
algebras of Coxeter groups for which such a homomorphism always exists.
\end{Example}

\begin{Remark}\label{not Mak}
As discussed in Remark~\ref{Mak}, Mak has shown that the cosets of the
reflection subgroups of $\G$ have a unique element of minimal length
with respect to the Bremke--Malle length function $\ell_0$ (see
Remark~\ref{length remark}). For each $\mu\in\SComp$ let $\E_\mu'$ be Mak's
set of distinguished coset representatives for $G_\mu$ and let
$E_\mu'=\sum_{e\in\E_\mu'}e\in R\G$. Define
$$\Sigma'(G_{r,n})=\sum_{\mu\in\SComp}RE_\mu'.$$
If $r>2$ then $\Sigma'(G_{r,n})$ is not, in general, a subalgebra of $R\G$. 
The smallest counter example occurs when $r=n=3$. 

Now suppose that $r=2$. Then $G_{2,n}$ is a Coxeter group of type $B_n$ and
Bonnaf\'e and Hohlweg~\cite{BonnafeHohlweg:hyper} have shown that
$\Sigma'(G_{2,n})$ is a subalgebra of~$RG_{2,n}$ and, moreover, that
$\Sigma'(G_{2,n})$ is isomorphic to
the Mantaci-Reutenauer algebra~\cite{MantaciReutenauer:descent}.  Now, the
algebras $\Sol(G_{2,n})$ and $\Sigma'(G_{2,n})$ are both free of rank
$2\cdot 3^{n-1}$, so it is natural to ask whether these algebras are
isomorphic. In fact, $\Sol(G_{2,n})\not\cong\Sigma'(G_{2,n})$ if $n>1$.
This can be proved by induction on~$n$ starting from the following
observation.  Bonnaf\'e and Hohlweg have shown 
in~\cite[Table V]{BonnafeHohlweg:hyper} that the following matrix is the
character table of the semisimple quotient of
$\Sigma'(G_{2,2})$.
$$\begin{array}{l|*{6}{c}}
  &(2) & (1^2) & (1,-1) & (-2) & (-1^2)\\\toprule
(2)   &1&&&&\\
(1^2) &1&2&&&\\
(1,-1)&1&2&2&&\\
(-2)  &1&.&.&2&&\\
(-1^2)&1&2&4&4&8&\\
\end{array}$$
Observe that the $\((-2),(-2)\)$--entry in this character table is
different to the corresponding entry in the character table of
$\Sol(G_{2,2})$ given in Example~\ref{character}. Therefore, $\Sol(G_{2,2})$
and $\Sigma'(G_{2,2})$ are not isomorphic algebras because they have
non--isomorphic maximal semisimple quotients.
\end{Remark}

We close this section by classifying the irreducible $\Sol_q(n)$--modules
over an arbitrary field. This classification is a direct generalization of
the corresponding results for the descent algebra of the symmetric
groups~\cite{AtkPfiWill} -- although our proofs are necessarily different
because there is no homomorphism from $\Sol_q(n)$ into the character
ring of~$\G$.

For $\lambda\in\Comp$ let
$N_{\Sym_n}(\Sym_\lambda)=\set{w\in\Sym_n|\Sym_\lambda=w^{-1}\Sym_\lambda
w}$ be the normalizer of $\Sym_\lambda$ in~$\Sym_n$.

\begin{Thm}\label{general irreducibles}
Suppose that $R$ is field, $q\in R$ and $\lambda\in\SPart$. Then
the following are equivalent:
\begin{enumerate}
	\item $d_{\lambda\lambda\lambda}(q)=0$;
  \item $q^{|\lambda|^-}[N_{\Sym_n}(\Sym_{\lamp}):\Sym_{\lamp}]=0$ in $R$;
	\item $E_\lambda\in\Rad\Sol_q(n)$;
  \item $E_\lambda$ is nilpotent; and,
	\item $I(\lambda)\cong I(\mu)$, for some $\mu\in\SPart$ with
    $|\Pi_\mu|>|\Pi_\lambda|$.
\end{enumerate}
\end{Thm}

\begin{proof}
By definition,
$$d_{\lambda\lambda\lambda}(q)
   =\sum_{\substack{d\in\D_{\lamp\lamp}\\\lambda=\lambda\cap d\lambda}}
        |T_{-\lambda\cap d(-\lambda)}|
   =\sum_{\substack{d\in\D_{\lamp\lamp}\\\lambda=\lambda\cap d\lambda}}
        q^{|\lambda|^-}
   = q^{|\lambda|^-}[N_{\Sym_n}(\Sym_{\lamp}):\Sym_{\lamp}],$$
since $|T_{-\lambda}|=q^{|\lambda|^-}$ and 
$T_{-\lambda\cap d(-\lambda)}=T_{-\lambda}$ if $\lambda=\lambda\cap d\lambda$. 
Hence, (a) and~(b) are equivalent. Further,~(c) and~(d) are equivalent by
Corollary~\ref{nilpotent}.

To complete the proof it is enough to show that
(a)$\implies$(c)$\implies$(e)$\implies$(a). In order to do this let
$\Sol_q(n)=\S_{\lambda_1}\supset\S_{\lambda_2}\supset\dots\supset
      \S_{\lambda_N}\supset0$ 
be the filtration of $\Sol_q(n)$ by two sided ideals which was constructed
in the proof of Theorem~\ref{irreducibles} using a total order~$>$
on~$\SComp$. Recall that $|\Pi_\mu|\ge|\Pi_\nu|$ whenever
$\mu>\nu$, for $\mu,\nu\in\SComp$. Then $\S_{\lambda_i}$ is a subalgebra
of $\Sol_q(n)$ which is also a quotient of $\Sol_q(n)$ since
$\S_{\lambda_i}\cong\Sol_q(n)/\S_{\lambda_{i+1}}$, for $1\le i\le N$.
Therefore, by Theorem~\ref{irreducibles}, every irreducible
$\S_{\lambda_i}$--module is isomorphic to $I(\mu)$ for some $\mu\in\SPart$
with $\mu\ge\lambda_i$, for $1\le i\le N$. In particular, every irreducible
$\S_\lambda$--module is isomorphic to $I(\mu)$ for some $\mu\ge\lambda$.

We can now return to the proof of the Theorem.

First, suppose (a) holds, so that $d_{\lambda\lambda\lambda}(q)=0$.  By
definition, if $\mu\in\SPart$ then $E_\lambda$ acts on $I(\mu)$ as
multiplication by $d_{\mu\lambda\mu}(q)$. By Lemma~\ref{equivalence}(a), if
$\mu>\lambda$ then $E_\lambda$ acts on $I(\mu)$ as multiplication by $0$,
whereas $E_\lambda$ acts on $I(\lambda)$ as multiplication by $0$ since
$d_{\lambda\lambda\lambda}(q)=0$. Therefore, $E_\lambda\in\Rad\S_\lambda$
and (c) holds because $\Rad\S_\lambda\subseteq\Rad\Sol_q(n)$.

Next, suppose that (c) holds. Then $E_\lambda$ belongs to the radical of
$\S_\lambda$. Now, $\S_\lambda\subset\S_{\lambda_{l-1}}$ so, as vector
spaces, $\Rad\S_\lambda=RE_\lambda+\Rad\S_{\lambda_{l-1}}$. On the other
hand, $\dim\S_\lambda=\dim\S_{\lambda_{l-1}}+1$, so it follows that the
algebras $\S_\lambda$ and $\S_{\lambda_{l-1}}$ have the same number of
irreducible modules. Hence, $I(\lambda)\cong I(\mu)$ for some signed partition
$\mu>\lambda$. That is, (e) holds.

Finally, assume that (e) holds. Then $I(\lambda)\cong I(\mu)$, for some
signed partition $\mu>\lambda$. Therefore, $E_\lambda$ acts on these
modules as multiplication by
$d_{\lambda\lambda\lambda}(q)=d_{\mu\lambda\mu}(q)$. Consequently,
$d_{\lambda\lambda\lambda}(q)=0$ by Lemma~\ref{equivalence}, so (a) holds.

This completes the proof of the Theorem.
\end{proof}

In the following Corollaries note that the integer
$d_{\lambda\lambda\lambda}(q)
    =q^{|\lambda|^-}[N_{\Sym_n}(\Sym_{\lamp}):\Sym_{\lamp}]$
is explicitly known by Theorem~\ref{general irreducibles} (and
Corollary~\ref{Emunu factorization}).

\begin{Cor}
Suppose that $R$ is a field and $q\in R$. Then
$$\set{I(\lambda)|\lambda\in\SPart\text{ and }
                 d_{\lambda\lambda\lambda}(q)\ne0}$$
is a complete set of pairwise non--isomorphic irreducible $\Sol_q(n)$--modules.
\end{Cor}

\begin{proof}
    This follows from Theorem~\ref{general irreducibles} and
    Theorem~\ref{irreducibles}.
\end{proof}

Similarly, combining the Theorem~\ref{general irreducibles} with
Corollary~\ref{nilpotent} and Corollary~\ref{radical}, we obtain the
general description of the radical of $\Sol_q(n)$ when $R$ is a field.

\begin{Cor}
Suppose that $R$ is a field and $q\in R$. Then
$$\set{E_\lambda-E_\mu|\lambda\in\SPart, \mu\in\SComp, \lambda\sim\mu
         \text{ and }\lambda\ne\mu}
\bigcup\set{E_\lambda|\lambda\in\SPart\text{ and }
            d_{\lambda\lambda\lambda}(q)=0}$$
is a basis of $\Rad\Sol_q(x)$.
\end{Cor}

Finally, we can use Theorem~\ref{general irreducibles} to describe the
radical and irreducible modules for each of the subalgebras of $\Sol_q(n)$
described in Proposition~\ref{subalgebras}. For brevity we state only the
following result.

\begin{Cor}
Suppose that $R$ is a field, $n\ge1$ and $q\in R$. Let $A$ be one of the
subalgebras $\Sol_q^+(n)$, $\Sol_q^\pm(n)$, $\Sol_q^1(n)$ of $\Sol_q(n)$.
Then $\Rad A=A\cap\Rad\Sol_q(n)$.
\end{Cor}

\section{The Hopf algebra of cyclotomic Solomon algebras}
In this section we fix $r>1$ and show that the direct sum of cyclotomic
Solomon algebras $\bigoplus_{n\ge0}\Sol(\G)$ is a concatenation Hopf
algebra, where $G_{r,0}=\{1_{G_{r,0}}\}$ is the trivial
group.  Further, this Hopf algebra is a Hopf subalgebra of the Hopf
algebra of colored permutations introduced by Baumann and
Hohlweg~\cite{BaumannHohlweg:coloredhopf}.

Most of the results in this section hold over an arbitrary integral domain,
however, the main results of this section (Theorem~\ref{Hopf} and
Corollary~\ref{subHopf}), hold only in characteristic zero. Consequently,
for this section we fix a field $\field$ of characteristic zero and we work
only over this field.  Thus, all tensor products are over $\field$, all
modules are $\field$-vector spaces and all algebras are $\field$--algebras.
In particular, the cyclotomic Solomon algebras $\Sol(\G)=\Sol_\field(\G)$
are $\field$--algebras.

We first recall some general facts about bialgebras and Hopf algebras. 

A \emph{$\field$--coalgebra} is a triple $(A,\delta,\epsilon)$ consisting
of a $\field$--vector space $A$ together with two linear maps $\delta\map
AA\otimes A$ (comultiplication) and $\epsilon\map A\field$ (the counit)
such that
$$(\delta\otimes \id_A)\circ\delta=(\id_A\otimes\delta)\circ\delta
\qquad\text{and}\qquad
(\epsilon\otimes \id_A)\circ\delta=(\id_A\otimes\epsilon)\circ\delta,$$
where $\id_A$ is the identity map on $A$.

A \emph{$\field$--bialgebra} is a coalgebra $(A,\delta,\epsilon)$ such
that $A$ is a $\field$--algebra and the structure maps $\delta\map
AA\otimes A$ and $\epsilon\map A\field$ are algebra homomorphisms. A
\emph{Hopf algebra} is a quadruple $(A,\delta,\epsilon,S)$ where
$(A,\delta,\epsilon)$ is a bialgebra and $S\map AA$ (the antipode) is a
linear map such that 
$\mu(S\otimes \id_A)\delta=\eta\epsilon=\mu(1\otimes S)\delta.$ 
Here $\mu\map{A\otimes A}A:(a,b)\mapsto ab$ is the multiplication map and
$\eta\map\field A;1\mapsto 1_A$ is the unit map for the algebra $A$.

Finally, a \emph{graded bialgebra} is a triple $(A,\delta,\epsilon)$ where
$A=\bigoplus_{n\in\N}A_n$ is $\N$--graded bialgebra and the maps $\delta$
and $\epsilon$ are graded (degree zero) vector space homomorphisms.  A
\emph{graded Hopf} algebra is a graded bialgebra which is equipped with an
antipode which is a graded vector space homomorphism of degree zero. A
graded bialgebra, or a graded Hopf algebra, $A=\bigoplus_{n\ge0}A_n$ is
\emph{connected} if $A_0=\field$.

Following Baumann and Hohlweg~\cite{BaumannHohlweg:coloredhopf}, we next
define the (graded connected) Hopf algebra of coloured permutations. This
will require some preparation. As a graded vector space this Hopf algebra
is the direct sum of the group algebras of groups~$\G$:
\[  \ColG := \bigoplus_{n\geq 0} \field G_{r,n}.\] 
We need some more notation before we can
describe the Hopf algebra structure on $\ColG$.

First, suppose that $m$ and $n$ are non--negative integers. Then
$G_{r,m}\times G_{r,n}$ is naturally isomorphic to the reflection subgroup
$G_{(m,n)}$ of~$G_{r,m+n}$. By identifying $G_{r,m}\times G_{r,n}$ and
$G_{(m,n)}$ we have an embedding 
$G_{r,m}\times G_{r,n}\hookrightarrow G_{r,m+n}$. Explicitly, this
embedding sends the generators $\{s_0,\dots,s_{m-1}\}$ of~$G_{r,m}$ to
$\{s_0,\dots,s_{m-1}\}$ in~$G_{r,m+n}$ and the generators
$\{s_0,\dots,s_{n-1}\}$ of $\G$ to $\{t_{m+1},s_{m+1},\dots,s_{m+n-1}\}$,
respectively.

By Proposition~\ref{row standard} there is a natural bijection between the
set $\E_{(m,n)}=\D_{(m,n)}$ of right coset representatives of $G_{(m,n)}$
in $\G$ and the set of row standard $(m,n)$--tableau. The product $*$
on the Hopf algebra $\ColG$ is the bilinear map determined by
$$u*v=\sum_{e\in\E_{(m,n)}}(u\times v)e=(u\times v)E_{(m,n)},$$
for $u\in G_{r,m}$, $v\in G_{r,n}$ and where $u\times v$ is multiplication
inside $G_{r,m+n}$ (the \textit{internal product} on $\ColG$). The product
$*$ on $\ColG$ is called the \emph{shuffle product}, or the \textit{external
product}, on~$\ColG$ because,  by
Proposition~\ref{row standard}, $\E_{(m,n)}$ is in bijection with the ways
of shuffling the two sets $\{1,\dots,m\}$ and $\{m+1,\dots,m+n\}$ together.
It is easy to check that $E_{(0)}=1_{G_{r,0}}\in\Sol(G_{r,0})$ is the unit
for the shuffle product.

To define the coproduct on $\ColG$ observe that for $m=0,\dots,n$ any element 
$g\in G_{r,n}$ can be written uniquely in the form 
$g=e_m^{-1}(g_{(m)}\times g_{(n)})$, where $g_{(m)}\in G_{r,m}$, 
$g_{(n)}\in G_{r,n}$ and $e_m\in\E_{(m,n)}$. Using this notation, the
\emph{coproduct} $\Delta$ on $\ColG$ is the linear map determined by
$$\Delta(g)=\sum_{m=0}^n g_{(m)}\otimes g_{(n)},$$
for $g\in\G$.

\begin{Example} In order to better
  distinguish between the elements $\G$ for different values of $n$ recall
  from the end of section~2 that there is a natural bijection between $\G$
  and the set of words 
  $\Words=\set{\underline\omega=\omega_1\dots\omega_n|
       \omega_i\in\rn\text{ and }\{|\omega_1|,\dots,|\omega_n|\}=\n}.$ 
  To give an example of the shuffle product and the coproduct on $\ColG$
  we identify $\G$ and $\Words$ using this bijection.

Suppose that $0\le a,b,c,d<r$. Then, using the identification above,
\begin{align*}
1\zeta^a\ 2 \zeta^b \ast 2\zeta^c\ 1\zeta^d 
  &=1\zeta^a2\zeta^b4\zeta^c3\zeta^d+1\zeta^a3\zeta^b4\zeta^c2\zeta^d   
   +1\zeta^a4\zeta^b3\zeta^c2\zeta^d\\
  &\phantom{=}+2\zeta^a3\zeta^b4\zeta^c1\zeta^d 
    +2\zeta^a4\zeta^b3\zeta^c1\zeta^d +3\zeta^a4\zeta^b2\zeta^c1\zeta^d\\
\intertext{and}
\Delta(2\zeta^a3\zeta^b1\zeta^c4\zeta^d) &= 
  \emptyset\otimes 2\zeta^a3\zeta^b1\zeta^c4\zeta^d
   +1\zeta^c \otimes1\zeta^a2\zeta^b3\zeta^d 
   + 2\zeta^a1\zeta^c \otimes 1\zeta^b2\zeta^d   \\ 
   &\phantom{=}+ 2\zeta^a3\zeta^b1\zeta^c\otimes 1\zeta^d  
   + 2\zeta^a3\zeta^b1\zeta^c4\zeta^d\otimes\emptyset,
\end{align*}
where $\emptyset$ is the empty word in $G_{r,0}$.
\end{Example}

As remarked above, $E_{(0)}=1_{G_{r,0}}$ is the multiplicative unit for the
shuffle product.  The counit of~$\ColG$ is the linear map
$\epsilon\map\ColG\field$ defined by 
\[ \epsilon(w) =\begin{cases} 1 & \text{ if } w=E_{(0)}\in G_{r,0} \\ 
                                 0 & \text{otherwise.} \end{cases}
\]

\begin{Thm}[%
   \protect{Baumann and Hohlweg~\cite[Theorem~1]{BaumannHohlweg:coloredhopf}}]
The triple $(\ColG,\Delta,\epsilon)$ is a graded connected bialgebra.
\end{Thm}

In fact, $(\ColG,\Delta,\epsilon)$ is a Hopf algebra at least when $\field$
is a field because every connected $\mathbb{N}$-graded $\field$-bialgebra
is a Hopf algebra; see \cite[Ex. 1, page 238]{Sweedler}.

We remind the reader that $r>1$ is fixed throughout this section.

\begin{Defn}
  The \emph{cyclotomic Hopf algebra} is the graded vector space
\[ \Sol(r)=\bigoplus_{n\geq 0} \Sol(G_{r,n}).\]
\end{Defn}

The cyclotomic Hopf algebra is naturally graded with $\Sol(r)_n=\Sol(\G)$
and, as a vector space, $\Sol(r)_n$ is finite dimensional with basis
$\set{E_\mu|\mu\in\SComp}$. For convenience, we set $E_n=E_{(n)}$, for
$n\in\Z$.

Our next aim is to show that $\Sol(r)$ is a Hopf subalgebra of $\ColG$. We
begin with a Lemma which generalizes~\ref{BBHT}(a).

\begin{Lemma}\label{coset fact}
  Suppose that $\alpha,\beta\in\SComp$ with $G_\alpha\subseteq G_\beta$.
  Then $\E^\beta_\alpha=\E_\alpha\cap G_\beta$ is a complete set of minimal
  length right coset representatives for $G_\alpha$ in~$G_\beta$ and
  $\E_\alpha=\E^\beta_\alpha\E_\beta$.
\end{Lemma}

\begin{proof} It is clear that $\E^\beta_\alpha$ is a complete set of right
coset representatives for $G_\alpha$ in~$G_\beta$. Moreover, by definition,
if $e\in\E^\beta_\alpha$ then $e$ is the unique element of minimal length
in the coset $G_\alpha e$. To prove the second statement observe that
$$\G=\coprod_{d\in\E_\beta} G_\beta d
  =\coprod_{d\in\E_\beta}\Big(\coprod_{e\in\E^\beta_\alpha} G_\alpha e\Big)d.$$
So, $\E^\beta_\alpha\E_\beta$ is a complete set of coset
representatives for $G_\alpha$ in~$\G$. Therefore, 
$\E_\alpha=\E^\beta_\alpha\E_\beta$ since the elements of both sides are of
minimal length in their respective cosets.
\end{proof}

\begin{Prop}
\label{Prop:shuffleprod}
Suppose that $\mu\in\Lambda_m^\pm$ and $\nu\in\Lambda_n^\pm$. Then 
\[E_{\mu}\ast E_{\nu} = E_{\mu\concat\nu}\in \Sol(G_{r,n+m})\]
where $\mu\concat\nu = (\mu_1, \ldots, \mu_{l}, \nu_1,\ldots, \nu_{k})$ is
the concatenation of two signed permutations. 
\end{Prop}

\begin{proof}
By definition, $E_\mu*E_\nu =(E_\mu\times E_\nu)E_{(m,n)}$
where, as above, we interpret $E_\mu\times E_\nu$ as an element of 
$\field G_{(m,n)}\subseteq\field\G$. Therefore, it is enough to prove that
$\E_{\mu\concat\nu}=\E^{(m,n)}_{\mu\times\nu}\E_{(m,n)}$. However, this
follows immediately from the previous Lemma because 
$G_{\mu\concat\nu}=G_\mu\times G_\nu\subseteq G_{(m,n)}$.
\end{proof}

Notice that the Proposition says that $\Sol(r)$ is a subalgebra of $\ColG$
and that, as an algebra, $\Sol(r)$ is freely generated by the elements
$\set{E_{\pm n}|n\ge 1}$.

\begin{Prop}\label{prop: coprod}
Suppose that $n$ is a positive integer. Then
\begin{enumerate}
\item $\Delta(E_n) = \sum\limits_{m=0}^n E_m\otimes E_{n-m}$;
\item $\Delta(E_{-n}) = \sum\limits_{m=0}^n E_{-m}\otimes E_{m-n}$.
\end{enumerate}
\end{Prop}

\begin{proof}
Part~(a) follows directly from the definitions. This result is well known
because $E_n=1_{\G}$ is the identity element of $\field\G$, so we omit
the details.  

For part (b), observe that $E_{-n}=F_{(n)}=\sum_{t\in T}t$. Therefore,
\begin{align*}
  \Delta(E_{-n})
  &=\sum_{\alpha=(\alpha_1,\dots,\alpha_n)\in\Z^r_n}
         \Delta(t_1^{\alpha_1}\dots t_n^{\alpha_n})\\
  &=\sum_{\alpha=(\alpha_1,\dots,\alpha_n)\in\Z^r_n}
    \sum_{m=0}^n t_1^{\alpha_1}\dots t_m^{\alpha_m}\otimes 
         t_1^{\alpha_{m+1}}\dots t_{n-m}^{\alpha_m}\\
  &=\sum_{m=0}^n \sum_{\substack{\beta\in\Z^r_m\\\gamma\in\Z^r_{n-m}}}
         t_1^{\beta_1}\dots t_m^{\beta_m}\otimes
         t_1^{\gamma_1}\dots t_{n-m}^{\gamma_{n-m}}\\
  &= \sum\limits_{m=0}^n E_{-m}\otimes E_{m-n},
\end{align*}
as required.
\end{proof}

We henceforth adopt the unusual convention that 
$\sum_{m=a}^bf(m)=\sum_{m=b}^af(m)$ if $b<a$. This allows us to write the
Proposition~\ref{prop: coprod} more compactly as
$\Delta(E_n)=\sum_{m=0}^{\mu_k} E_m\otimes E_{\mu_k-m}$, for $n\in\Z$.

As the coproduct is an algebra homomorphism $\ColG\to\ColG\otimes\ColG$ it
follows from the last two Propositions that $\Sol(r)$ is a sub-bialgebra of
$\ColG$.

Let $\P$ be a set of non-commuting indeterminates over~$\field$. The
\emph{concatenation Hopf algebra} on $\P$ is the free associative
$\field$-algebra $\field\<\P\>$ on~$\P$ with 
counit $\epsilon$, where $\epsilon(f(\P))=f(0)$ is the constant term of
$f(\P)\in\field\<\P\>$, coproduct $\delta(p)=p\otimes1+1\otimes p$ for any
$p\in \P$, and antipode $S$ determined by $S(p_1\dots p_k)=(-1)^kp_k\dots
p_1$, for $p_1,\dots,p_k\in \P$. Any function $\deg\map\P\N$ extends to a
degree function on the monomials in $\field\<\P\>$ by setting
$\deg(p_1\dots p_k)=\deg(p_1)+\dots+\deg(p_k)$. In this way, 
$\field\<\P\>=\bigoplus_{n\ge0}\field\<\P\>_n$
becomes a graded connected Hopf algebra, where $\field\<\P\>_n$ is the space
of homogeneous polynomials $p_1\dots p_k$ in~$\P$ with 
$\deg(p_1\dots p_k)=n$.

We can now prove the main result of this section. Up until now we have not
used the assumption that $\field$ is a field of characteristic zero. This
assumption is necessary, however, for the proof of the following Theorem.

\begin{Thm}\label{Hopf}
  Suppose that $\field$ is a field of characteristic zero. Then 
  $(\Sol(r),\Delta,\epsilon)$ is isomorphic to the graded
  connected concatenation Hopf algebra $\field \<\P\>$ on a set of
  non-commuting indeterminates $\P=\set{P_n|n\in\Z\setminus\{0\}}$ 
  where $\deg P_{\pm n}=n$, for $n>0$.
\end{Thm}

\begin{proof} 
Our argument is modeled on the proof of  \cite[Theorem~2.1]{MalvenutoReutenauer:descent}. 

Let $x$ be a formal variable and consider the algebra $\Sol(r)\llbracket
x\rrbracket$ of formal power series in $x$ over $\Sol(r)$, where $x$
commutes with $\Sol(r)$. For each positive integer~$n$ define elements
$P_{\pm n}\in\Sol(r)\llbracket x\rrbracket$ using the generating series
\begin{align*} \label{eq:defposp}
  \sum_{n>0} P_n x^n &= \log(1+ E_1 x + E_2 x^2 + \cdots)\\
\intertext{and}
\sum_{n>0} P_{-n} x^n &= \log(1+ E_{-1} x + E_{-2} x^2 + \cdots).
\end{align*}
A straightforward calculation using Proposition~\ref{Prop:shuffleprod} and 
the Taylor series expansion of $\log(1+t)$ shows that
$$P_n=\sum_{\alpha\in\Comp}\frac{(-1)^{\ell(\alpha)-1}}{\ell(\alpha)}E_\alpha
\quad\text{and}\quad
P_{-n}=\sum_{-\alpha\in\Comp}\frac{(-1)^{\ell(\alpha)-1}}{\ell(\alpha)}E_\alpha.$$
(Recall that $\ell(\alpha)$ is the number of non-zero parts in $\alpha$.)
Therefore, $P_{n}, P_{-n}\in\Sol(\G)$ are homogeneous of degree~$n$; in
particular, $P_{\pm n}\in\Sol(r)$, for all $n>0$. Consequently, the elements 
$\set{P_{\pm n}|n>0}$ generate a subalgebra of $\Sol(r)$.

Similarly, since
$\sum_{n\ge0}E_{\pm n}x^n=\exp(\sum_{n>0}P_{\pm n}x^n)$, another
completely formal calculation using the Taylor series expansion of $\exp(x)$ and
Proposition~\ref{Prop:shuffleprod} shows that if $n>0$ then
$$E_n=\sum_{\alpha\in\Comp}\frac1{\ell(\alpha)!}P_\alpha
\quad\text{and}\quad
E_{-n}=\sum_{-\alpha\in\Comp}\frac1{\ell(\alpha)!}P_\alpha,$$
where we set $P_\alpha=P_{\alpha_1}*\dots*P_{\alpha_k}$, for 
$\alpha=(\alpha_1,\dots,\alpha_k)\in\pm\Comp$. Therefore, by the last
paragraph, the set
$\P=\set{P_n|n\in\Z\setminus\{0\}}$ freely generates 
$\Sol(r)$ as an algebra. That is, $\Sol(r)=\<P_{\pm n}\mid n>0\>$ as an algebra.

We claim that
$\Delta(P_n)=P_n\otimes1+1\otimes P_n$, for $n\in\Z\setminus\{0\}$. This
will complete the proof because it shows that these elements generate a
concatenation Hopf algebra $\field\<\P\>$ inside $\Sol(r)$. Starting from
the definition of $P_{\pm n}$ we have that 
$$\sum_{n>0}\Delta(P_{\pm n}x^n)
    =\Delta\(\sum_{n\ge0}\log(\sum_{n\ge0}E_{\pm n}x^n)\)
    =\log\(\sum_{n\ge0}\Delta(E_{\pm n})x^n\),
$$
where the last equality follows by the linearity of Taylor expansions since
$\Delta$ is an algebra homomorphism. Using
Proposition~\ref{prop: coprod} to expand the right hand side of the last
equation, exactly as in the proof
of \cite[(2.9)]{MalvenutoReutenauer:descent}, shows that
$\Delta(P_{\pm n})=P_n\otimes1+1\otimes P_n$. This proves our claim and so 
completes the proof. 
\end{proof}

\begin{Cor}\label{subHopf}
Suppose that $r>1$. Then the graded vector space $\Sol(r)$ equipped with
the product $\ast$, coproduct $\Delta$, unit $E_0$ and
counit $\epsilon$, is a graded connected Hopf subalgebra of $\ColG$. 
\end{Cor}

\section{A second bialgebra structure on $\Sol(r)$}
In this section we show that the cyclotomic Hopf algebra $\Sol(r)$ has a
second bialgebra structure with the same coproduct $\Delta$ as in
section~9, but where the product is inherited from multiplication in
the groups $\G$, for $r,n\ge0$. More precisely, the \textit{internal product}
is the unique bilinear map $\cdot\map\ColG\ColG$ such that if
$w\in G_{r,m}$ and $v\in\G$ then 
$$ w\cdot v=\begin{cases}
        wv,&\text{if } n=m,\\
        0,&\text{otherwise}.
\end{cases}$$
We frequently abuse notation and write $xy=x\cdot y$, for
$x,y\in\ColG$. 

As each of the group algebras $\field\G$ are associative algebras it
follows that $(\ColG,\cdot)$ is an associative algebra.  Note, however, that
$(\ColG,\cdot)$ does not have a multiplicative unit, so we cannot expect to
obtain a second Hopf algebra structure on $\Sol(r)$ in this way. Note also
that the internal product~$\cdot$ does not respect the grading on
$\ColG=\bigoplus_n\field\G$. 

By Theorem~\ref{cyc Solomon}, $\Sol(r)$ is a subalgebra of the algebra
$(\ColG,\cdot)$. One can easily check that $\epsilon$ is an algebra homomorphism 
on $(\ColG, \cdot)$, whereas $\Delta$ is not an algebra homomorphism on $\ColG$, see \cite[Remark~5.15]{Malvenuto}. However, we will show that $\(\Sol(r),\cdot,\Delta\)$ is a
bialgebra. To prove this we need only show that $\Delta$ is
an algebra homomorphism with respect to the internal product. The argument
that we give generalizes that used by Malvenuto
\cite[Remark~5.15]{Malvenuto} to prove the analogous statement for the
descent algebra of the symmetric group. We start with some new definitions.

A \textbf{pseudo signed composition} of $n$ is an element
$\c=(c_1,c_2,\ldots,c_k)\in\Z^k$, for some $k>0$, such that
$|\c|=|c_1|+|c_2|+\cdots |c_k|=n$. A \textbf{pseudo composition} is an
element of $\N^k$, for some $k>0$. The difference between (signed)
compositions and pseudo (signed) compositions is that pseudo (signed) 
compositions can contain zeros. If $\c$ is a pseudo signed composition let
$\bar\c$ be the signed composition obtained by omitting the zeros from
$\c$. For example, if $\c=(-2,0,3,0,1)$ then $\bar\c=(-2,3,1)$.

If $\c\in\Z^k$ is a pseudo signed composition then set
$E_\c=E_{\bar\c}$. If $\c,\c'\in\Z^k$ are two pseudo signed composition of
the same length then $\c+\c'\in\Z^k$, where addition is defined
componentwise.  We extend the operation of concatenation to pseudo
signed compositions in the obvious way so that if $\c\in\Z^k$ and
$\c'\in\Z^l$ then $\c\concat\c'\in\Z^{k+l}$.

Two integers $c$ and $c'$ are \textbf{sign equivalent}, and we write
$c\sequiv c'$, if $cc'>0$.
Similarly, two (pseudo) signed compositions $\c=(c_1,\dots,c_k)$ and
$\c'=(c'_1,\dots,c'_k)$ are sign equivalent if $c_i\sequiv c'_i$, for
$i=1,\dots,k$. Again, we write $\c\sequiv\c'$.

\begin{Prop}\label{delta}
Suppose that $\mu\in\SComp$ and that $\len(\mu)=k$. Then
$$\Delta(E_\mu)=\sum_{\substack{\c'\sequiv\c''\in\Z^k\\\mu=\c'+\c''}}
              E_{\c'}\otimes E_{\c''}.$$
\end{Prop}

\begin{proof}
  We argue by induction on~$k$. As $\Delta(E_0)=E_0\otimes E_0$ the
  case $k=0$ is clear. So we may assume that $k>0$. Let
  $\nu=(\mu_1,\dots,\mu_{k-1})$ so that $\mu=\nu\concat(\mu_k)$. Then, by
  Proposition~\ref{Prop:shuffleprod} and Proposition~\ref{prop: coprod},
  \begin{align*}
    \Delta(E_\mu)&=\Delta(E_\nu* E_{\mu_k})
                  =\Delta(E_\nu)*\Delta(E_{\mu_k})\\
                  &=\Big(\sum_{\substack{\c'\sequiv\c''\in\Z^{k-1}\\\nu=\c'+\c''}}
                        E_{\c'}\otimes E_{\c''}\Big)*
                        \Big(\sum_{m=0}^{\mu_k} E_m\otimes
                        E_{\mu_k-m}\Big),
  \intertext{by induction on $k$. (If $\mu_k<0$ then recall our unusual convention  
  for summations from after Proposition~\ref{prop: coprod}.) Therefore, using 
  Proposition~\ref{Prop:shuffleprod} for the second equality,}
  \Delta(E_\mu)&=\sum_{\substack{\c'\sequiv\c''\in\Z^{k-1}\\\nu=\c'+\c''}}
          \sum_{m=0}^{\mu_k} E_{\c'}*E_m\otimes E_{\c''}*E_{\mu_k-m}\\
          &=\sum_{\substack{\c'\sequiv\c''\in\Z^{k-1}\\\nu=\c'+\c''}}
          \sum_{m=0}^{\mu_k} E_{\c'\concat(m)}\otimes E_{\c''\concat(\mu_k-m)}\\
          &=\sum_{\substack{\c'\sequiv\c''\in\Z^k\\\mu=\c'+\c''}}
              E_{\c'}\otimes E_{\c''}
  \end{align*}
  as required.
\end{proof}

Let $k,l>0$ be positive integers and let $\Mat_{kl}(\Z)$ be the set of
$k\times l$ integer matrices. If $M\in\Mat_{kl}(\Z)$ let
$\row(M)=(r_1,\dots,r_k)$ be the pseudo composition where $r_i$ is the sum
of the absolute values of the entries in row $i$ of $M$, for $1\le i\le k$. Similarly, let
$\col(M)=(c_1,\dots,c_k)$ be the pseudo composition where $c_j$ is the sum
of the absolute values of the entries in column~$j$ of $M$.  Finally, if
$M\in\Mat_{kl}(\Z)$ let $\comp(M)$ be the signed composition obtained by
listing the non--zero entries in $M$ in order, from left to right and then
top to bottom; thus, if $M=(m_{ij})$ then
$\comp(M)=\bar{(m_{11},\dots,m_{1l},m_{21},\dots,m_{k1},\dots,m_{kl})}$.

If $\c=(c_1,\dots,c_k)$ is a pseudo signed composition then define
$\c^+=(|c_1|,\dots,|c_k|)$. In the next definition we are most interested
in the case when $\mu$ and $\nu$ are signed compositions. We include pseudo
signed compositions in the definition because they are needed in the proof
of Theorem~\ref{bialgebraII} below.

\begin{Defn}\label{matrices}
Suppose that $\mu=(\mu_1,\dots,\mu_k)$ and $\nu=(\nu_1,\dots,\nu_l)$ are
pseudo signed compositions of~$n$. Let 
$$\M_{\mu\nu}=\SET[50]M=(m_{ij})\in\Mat_{kl}(\Z)|
               $\row(M)=\mu^+$, and $\col(M)=\nu^+$,
             $m_{ij}\le 0$ if $\mu_i<0$ or if $\nu_j<0$,\\
             and $m_{ij}\geq 0$ if $\mu_i>$ and $\nu_j>0$
               |.$$
Suppose now that $M=(m_{ij})\in\M_{\mu\nu}$. The \textbf{weight} of $M$
is the non-negative integer 
$$\wt(M)=-\sum_{\substack{i: \mu_i<0\\j : \nu_j<0}} m_{ij},$$ 
where in the sum $1\le i\le k$ and $1\le j\le l$ (note that $m_{ij}\le 0$
for all such $i,j$). If $\mu$ and $\nu$ are signed compositions let $\T_M$
be the unique row semistandard tableau in~$\rssTab(\mu,\nu)$ such that $j$
appears $|m_{ij}|$ times in row~$i$ of~$\T$, for $1\le i\le k$ and 
$1\le j\le l$. 
\end{Defn}

Note that if~$\mu$ and~$\nu$ are compositions and $M=(m_{ij})\in\M_{\mu\nu}$
then $\wt(M)=0$ and $m_{ij}\ge0$, for $1\le i\le\len(\mu)$ and $1\le
j\le\len(\nu)$.

\begin{Prop}\label{combII}
  Suppose that $\mu$ and $\nu$ are signed compositions of $n$. Then
 \[E_\mu E_\nu = \sum_{M\in\M_{\mu\nu}} r^{\wt(M)} E_{\comp(M)}\]
 \end{Prop}

\begin{proof}By Theorem~\ref{cyc Solomon}, 
  $E_\mu E_\nu=\sum_{d\in\D_{\mup\nup}}|T_{-\mu\cap d(-\nu)}|\,E_{\mu\cap d\nu}.$
  Therefore, to prove the Proposition it is enough
  to show that there exists a bijection 
  $\M_{\mu\nu}\to\D_{\mup\nup};M\mapsto d_M$ such that
  $\comp(M)=\mu\cap d_M\nu$ and $r^{\wt(M)}=|T_{-\mu\cap d_M(-\nu)}|$.

  First, observe that the map $\M_{\mu\nu}\to\rssTab(\mu,\nu);M\mapsto\T_M$
  is a bijection because its inverse is the map which sends the tableau
  $\T\in\rssTab(\mu,\nu)$ to $M_\T=(m_{ij})$, where $|m_{ij}|$ is the
  number of times that~$j$ appears in row~$i$ of~$\T$, and where the sign
  of~$m_{ij}$ is determined by the constraints on $\M_{\mu\nu}$. Next,
  by~(\ref{double Sym}), the map 
  $\rssTab(\mu,\nu)\to\D_{\mup\nup};\T\mapsto d_{\T^*}$ is a bijection.
  Hence, the map 
  $$\M_{\mu\nu}\to\D_{\mup\nup}; M\mapsto d_M=d_{\T_M^*}$$
  is a bijection. 

  Fix $M\in\M_{\mu\nu}$. Then $\wt(M)=\wt(\T_M)$ in the notation of
  Corollary~\ref{Emunu factorization}, so that
  $r^{\wt(M)}=|T_{-\mu\cap d_M(-\nu)}|$. Hence, it remains to prove that
  $\comp(M)=\mu\cap d_M\nu$. The permutation $d_M$ is determined by the row
  semistandard tableau $\T$ which, by the last paragraph, also
  determines~$M=(m_{ij})$. If $m_{ij}\ne0$ then $|m_{ij}|$ is equal to the
  number of times that~$j$ appears in row~$i$ of~$\T$. Writing
  $G_\mu=G_{\mu_1}\times\dots\times G_{\mu_k}$ and
  $G_\nu=G_{\nu_1}\times\dots\times G_{\nu_l}$, and abusing notation
  slightly, we see that $m_{ij}$ computes the intersection of $G_{\mu_i}$
  with $d_MG_{\nu_j}d_M^{-1}$; more precisely,
  $$G_{\mu_i}\cap d_MG_{\nu_j}d_M^{-1}
       \cong\begin{cases}
         G(r,1,m_{ij}),&\text{if }m_{ij}>0,\\
         \Sym_{-m_{ij}},&\text{if }m_{ij}<0.
  \end{cases}$$
  Comparing this with the recipe given in the proof of 
  Lemma~\ref{mud cap nu} for computing $\mu\cap d_M\nu$ we see that 
  $\comp(M)=\mu\cap d_M\nu$, as required.
\end{proof}

Garsia--Remmel~\cite[Prop. 1.1]{GarsiaReutenauer} (see also
\cite[\Sect4]{GarsiaRemmel}), proved the analogue of this result for the
Solomon algebras of the symmetric groups. This is equivalent to the special
case of Proposition~\ref{combII} when~$\mu$ and~$\nu$ are both compositions
of~$n$. If $\mu,\nu\in\Comp$ then the bijection
$\M_{\mu\nu}\bijection\D_{\mu\nu}$ is well--known; see, for example,
\cite[Theorem~1.3.10]{JamesKerber}.

\begin{Example}\label{multiplication example 2}
   As in Example~\ref{row semistandard}, suppose that $\mu=(3,-2)$ and
   $\nu=(-2^2,1)$. The following table lists all of the elements of
   $\M_{\mu\nu}$, together with the associated signed composition and
   row semistandard $\mu$--tableau of type~$\nu$ and the weight of the
   matrix.
   $$\begin{array}{lllc}
     \multicolumn1c M & \comp(M)  & \multicolumn1c{\T_M} & \wt(M) \\\toprule
     \mat{-2&-1&0\\0&-1&-1} & (-2,-1^3) & \Tab(1&1&2\cr2&3) & 1 \\[10pt]
     \mat{-2&0&1\\0&-2&0} & (-2,1,-2) & \Tab(1&1&3\cr2&2)& 2 \\[10pt]
     \mat{-1&-2&0\\-1&0&-1} & (-1,-2,-1^2) & \Tab(1&2&2\cr1&3)& 1 \\[10pt]
     \mat{-1&-1&1\\-1&-1&0} & (-1^2,1,-1^2)& \Tab(1&2&3\cr1&2)& 2 \\[10pt]
     \mat{0&-2&1\\-2&0&0} & (-2,1,-2) & \Tab(2&2&3\cr1&1) & 2\\
\end{array}$$
The reader might like to compare this table with the one given in
Example~\ref{row semistandard}.

Combining the information above with Proposition~\ref{combII} shows that 
$$E_{(3,-2)}E_{(-2^2,1)}
  =2r^2E_{(-2,1,-2)}+rE_{(-2,-1^3)}+rE_{(-1,-2,-1^2)}+r^2E_{(-1^2,1,-1^2)}.$$
This calculation agrees with Example~\ref{multiplication example 1}, as it
must.
\end{Example}

Suppose that $M'=(m'_{ij}), M''=(m''_{ij})\in\Mat_{kl}(\Z)$, for some $k,l>0$. Then $M'$ and $M''$ are \textbf{signed equivalent}, and we write $M'\sequiv M''$, if $m'_{ij}\sequiv m''_{ij}$, for $1\le i\le k$ and $1\le j\le l$.

We can now prove the main result of this section. 

\begin{Thm}\label{bialgebraII}
  Suppose that $r>1$. Then $\Sol(r)$ equipped with product $\cdot$,
  coproduct $\Delta$ and counit $\varepsilon$, is a bialgebra. 
\end{Thm}

\begin{proof}As remarked at the beginning of this section, it remains to show that
  the coproduct $\Delta\map{\Sol(r)}{\Sol(r)\otimes\Sol(r)}$ is an algebra 
  homomorphism with respect to the internal product. By linearity it
  is enough to show that 
  $$\Delta(E_\mu E_\nu) = \Delta(E_\mu)\Delta(E_\nu),$$ 
  for all signed compositions $\mu$ and $\nu$. Further, we may assume that 
  $|\mu|=|\nu|$ since otherwise both sides of this equation are zero. Let 
  $k=\len(\mu)$ and $l=\len(\nu)$ and for $M\in\M_{\mu\nu}$ let
  $\len(M)=\len(\comp(M))$. Then, by Proposition~\ref{combII} and
  Proposition~\ref{delta},
\begin{align*}
\Delta(E_\mu E_\nu) &= \sum_{M\in\M_{\mu\nu}} r^{\wt(M)}\Delta(E_{\comp(M)})\\
  &=\sum_{M\in\M_{\mu\nu}}
    \sum_{\substack{\c'\sequiv\c''\in\Z^{\len(M)}\\\comp(M)=\c'+\c''}}
        r^{\wt(M)} E_{\c'}\otimes E_{\c''},
\end{align*}
For the moment, fix a matrix $M\in\M_{\mu\nu}$ and $\c',\c''\in\Z^{\len(M)}$ 
such that $\c'\sequiv\c''$ and $\comp(M)=\c'+\c''$. Since $\c'\sequiv\c''$
there exist unique matrices $M'=(m_{ij}'),M''=(m_{ij}'')\in\Mat_{kl}(\Z)$ such
that $M=M'+M''$, $M\sequiv M'\sequiv M''$, $\comp(M')=\bar{\c'}$ and
$\comp(M'')=\bar{\c''}$. Note that $\wt(M)=\wt(M')+\wt(M'')$ since $M'\sequiv M''$. 
Therefore, the last equation becomes
$$\Delta(E_\mu E_\nu) 
    =\sum_{M\in\M_{\mu\nu}} \sum_{\substack{M'\sequiv M''\\M'+M''=M}} 
       r^{\wt(M')} E_{\comp(M')}\otimes r^{\wt(M'')} E_{\comp(M'')}
$$
For each pair $M'$ and $M''$ in the second sum let $\mu'=\row(M')$ and
$\mu''=\row(M'')$. Then~$\mu'$ and $\mu''$ are 
\textit{pseudo} signed compositions such that $\mu=\mu'+\mu''$
and $\mu'\sequiv\mu''$. Similarly, $\nu'=\col(M')$ and $\nu''=\col(M'')$ are
\textit{pseudo} signed compositions such that $\nu=\nu'+\nu''$ and
$\nu'\sequiv\nu''$. By signed equivalence, $M'\in\M_{\mu'\nu'}$ and
$M''\in\M_{\mu''\nu''}$. Moreover, $M'$ and $M''$ run through
$\M_{\mu'\nu'}$ and $\M_{\mu''\nu''}$, respectively, for all possible
$\mu',\mu'',\nu'$ and $\nu''$, as $M$ runs through $\M_{\mu\nu}$. Observe
that if $M'\in\M_{\mu'\nu'}$ and $M''\in\M_{\mu''\nu''}$, 
for~$\mu',\mu'',\nu'$ and~$\nu''$ as above, then $M'\sequiv M''$ since
$\mu'\sequiv\mu''$ and $\nu'\sequiv\nu''$.  Therefore, we can 
reverse the order of summation in the last displayed equation to obtain
\begin{align*}
\Delta(E_\mu E_\nu) 
  &=\sum_{\substack{\mu'\sequiv\mu''\\\mu=\mu'+\mu''\\\nu'\sequiv\nu''\\\nu=\nu'+\nu''}}
    \sum_{\substack{M'\in\M_{\mu'\nu'}\\M''\in\M_{\mu''\nu''}}}
    r^{\wt(M')}E_{\comp(M')}\otimes r^{\wt(M'')}E_{\comp(M'')}\\
  &=\sum_{\substack{\mu'\sequiv\mu''\\\mu=\mu'+\mu''\\\nu'\sequiv\nu''\\\nu=\nu'+\nu''}}\!\!\!\!
    \Big(\sum_{M'\in\M_{\mu'\nu'}}\!\!r^{\wt(M')}E_{\comp(M')}\Big)
    {\otimes}
    \Big(\sum_{M''\in\M_{\mu''\nu''}}\!\!r^{\wt(M'')}E_{\comp(M'')}\Big)\\
  &=\Big(\sum_{\substack{\mu'\sequiv\mu''\\\mu=\mu'+\mu''}}
           E_{\mu'}\otimes E_{\mu''}\Big)
    \Big(\sum_{\substack{\\\nu'\sequiv\nu''\\\nu=\nu'+\nu''}}
           E_{\nu'}\otimes E_{\nu''}\Big)\\
  &=\Delta(E_\mu)\Delta(E_\nu),
\end{align*}
where the last two equalities follow by Proposition~\ref{combII} and
Proposition~\ref{delta} respectively. This completes the proof.
\end{proof}

\section*{Acknowledgments}
Many of the results in this paper were inspired by extensive computer
calculations using programs written using \textsc{Gap}~\cite{GAP}. We thank
N.~Bergeron and M.~Aguiar for useful conversations and the referee for
their careful reading of our manuscript.



\end{document}